\documentclass[notitlepage,10pt]{article}
\usepackage{amssymb,amsbsy,amsmath,mathrsfs}
\usepackage{pgfplots}
\usepackage{comment}
\pgfplotsset{compat=newest}

\catcode`\@=11
\@addtoreset{equation}{section}

\catcode`\@=12

\newtheorem{Theorem}{Theorem}[section]
\newtheorem{Lemma}[Theorem]{Lemma}

\newtheorem{Corollary}[Theorem]{Corollary}

\def\Proof{\noindent{{\it Proof. }}}
\def\QED{\hfill {$\square$}\goodbreak\medskip}

\def\N{{\mathbb N}}
\def\R{{\mathbb R}}

\def\Z{{\mathbb Z}}
\def\E{{\mathscr E}}
\def\F{{\mathcal F}}

\newcommand{\de}{\,\mathrm{d}}

\newcommand{\n}{\mathbf{n}}
\newcommand{\K}{\mathcal{K}}

\newcommand{\Q}{\mathcal{Q}}
\newcommand{\A}{\mathcal{A}}
\newcommand{\Le}{\mathcal{L}}

\title{On the dynamics of a charged particle in\\ magnetic fields with cylindrical symmetry}

\author{Paolo Caldiroli\footnote{Dipartimento di Matematica, Universit\`a di Torino, via Carlo Alberto, 10 -- 10123 Torino, Italy. Email: \tt{paolo.caldiroli@unito.it}}, Gabriele Cora\footnote{Dipartimento di Matematica, Universit\`a di Torino, via Carlo Alberto, 10 -- 10123 Torino, Italy. Email: \tt{gabriele.cora@unito.it}}}

\date{}

\begin{document}
\maketitle

\begin{abstract}
\noindent
We study the motion of a charged particle under the action of a magnetic field with cylindrical symmetry. In particular we consider magnetic fields with constant direction and with magnitude depending on the distance $r$ from the symmetry axis of the form $1+Ar^{-\gamma}$ as $r\to\infty$, with $A\ne 0$ and $\gamma>1$. With perturbative-variational techniques, we can prove the existence of infinitely many trajectories whose projection on a plane orthogonal to the direction of the field describe bounded curves given by the superposition of two motions: a rotation with constant angular speed at a unit distance about a point which moves along a circumference of large radius $\rho$ with a slow angular speed $\varepsilon$. The values $\rho$ and $\varepsilon$ are suitably related to each other. This problem has some interest also in the context of planar curves with prescribed curvature. 
\smallskip

\smallskip

\end{abstract}

\section{Introduction and statement of the result}

This paper concerns the dynamics of a particle in a magnetic field. This subject is quite important in several areas of physics, such as condensed matter theory \cite{ForZim}, accelerator physics \cite{Mic}, magnetobiology \cite{HaaBroVenTho} and plasma physics \cite{Pir}. However, despite the relevance and the variety of contexts, it seems that for this problem, in concrete magnetic fields, only few rigorous results are available in the literature (see \cite{AguGinPer, AguLuqPer, CalGui, LuqPer, Tru}).

From classical physics, in presence of an external magnetic field ${\bf{B}}$, the motion $t\mapsto q(t)\in\R^{3}$ of a particle of mass $m$ and charge $e$ is driven by the Lorentz force, according to the equation
\begin{equation}
\label{eq:Lorentz}
m\ddot q=e\dot q\wedge \bf{B}.
\end{equation}
When $\bf{B}$ is a constant field, namely ${\bf{B}}=B\nu$ for some versor $\nu$ and $B$ positive constant, the system is integrable and trajectories are helices with axis a line of the field and Larmor ray given by $r_{L}=\frac{mv_{\bot}}{|e|B}$, where $v_{\bot}$ is the speed orthogonal to the field direction. More precisely, the projection $q_{\perp}(t)$ of $q(t)$ on a plane orthogonal to $\nu$ moves on a circle of radius $r_{L}$ with constant angular speed $|e|B/m$.

The situation can drastically change if one considers a variable magnetic field, even in the case of constant direction and magnitude close to a constant. In this work we consider a class of magnetic fields with cylindrical symmetry, namely we assume that $\bf{B}$ has a constant direction, which can be taken as the third axis, and the magnitude of $\bf{B}$ depends just on the distance from the third axis. This means that 
\begin{equation}
\label{eq:B3}
{\bf{B}}(q_{1},q_{2},q_{3})=(0,0,B(q_{1},q_{2}))\quad\forall (q_{1},q_{2},q_{3})\in\R^{3}
\end{equation}
with $B$ depending just on $\sqrt{q_{1}^{2}+q_{2}^{2}}$. For $\bf{B}$ of the form (\ref{eq:B3}), problem (\ref{eq:Lorentz}) reduces to the system
\begin{equation}
\label{eq:Lorentz2}
\left\{\begin{array}{l}\ddot q_{1}=\frac{e}{m}B(q_{1},q_{2})\dot q_{2}\\
\ddot q_{2}=-\frac{e}{m}B(q_{1},q_{2})\dot q_{1}\\
\ddot q_{3}=0~\!.
\end{array}\right.
\end{equation}
In this case, again the component of $q(t)$ in the direction of $\bf{B}$, i.e. $q_{3}(t)$, describes a uniform motion, but the projection of $q(t)$ on the horizontal plane follows a trajectory which in general can be far from being circular or closed or also globally bounded, even when $B$ is close to a constant. 

Considering that the study of the dynamics of a single charged particle in a magnetic field constitutes a first elementary model for the motion of plasma, the search for trapped trajectories (in the directions orthogonal to the magnetic field) seems properly motivated. For this reason we are interested in motions $q(t)$ whose projection on the horizontal plane, given by $t\mapsto (q_{1}(t),q_{2}(t))=:\tilde{q}(t)\in\R^{2}$, draws a globally bounded, possibly closed, planar curve. We will say that such motions are \emph{bounded}. When the projection $\tilde{q}(t)$ defines a closed curve, we will say that the corresponding motion is \emph{periodic}.

With respect to this goal, we can reduce ourselves to study the planar system defined by the first two equations of (\ref{eq:Lorentz2}). Setting $v(t)=\tilde{q}(-mt/e)$, the planar system for $\tilde{q}(t)$ becomes
\begin{equation}
\label{eq:Lorentz3}
\left\{\begin{array}{l}v\colon\R\to\R^{2}\\
\ddot v=iB(v)\dot v\end{array}\right.
\end{equation}
where we identify $\R^{2}$ with the complex plane; the multiplication by the imaginary unit yields a counter-clockwise rotation of $\pi/2$.

Actually, we can change the system in (\ref{eq:Lorentz3}) into an equivalent \emph{scalar} problem of geometrical type, concerning curves with prescribed mean curvature. To this aim, recall that for a regular curve $u\colon\R\to\R^{2}$ of class $C^{2}$, the curvature at $u(t)$ is given by
\begin{equation}
\label{mathcal-K}
\mathcal{K}(u(t))=\frac{i\dot u(t)\cdot\ddot{u}(t)}{|\dot u(t)|^{3}}~\!.
\end{equation}
If $v$ solves (\ref{eq:Lorentz3}) then $|\dot v|$ is constant (this can be obtained just multiplying the equation by $\dot v$, that implies $\frac{\de}{\de t}|\dot v|^{2}=0$). In particular, if a solution $v$ of (\ref{eq:Lorentz3}) verifies also $|\dot v|\equiv 1$, then $B(v)$ equals the curvature of the curve parametrized by $v$. On the other hand, if  $u\colon\R\to\R^{2}$ is a bounded solution to
\begin{equation}
\label{eq:curvature}
\frac{i\dot u\cdot\ddot{u}}{|\dot u|^{3}}=B(u)
\end{equation}
then one can find an increasing diffeomorphism $g$ of $\R$ onto $\R$ such that $u\circ g$ solves (\ref{eq:Lorentz3}) (see Lemma \ref{L:curvature}). In this way, for a vector field $\mathbf{B}$ with cylindrical symmetry, i.e., of the form (\ref{eq:B3}), we can covert the problem of bounded trajectories of (\ref{eq:Lorentz}) into the geometrical problem about bounded curves with prescribed curvature $B=|\mathbf{B}|$. 

There are many contributions to the issue of planar $K$-loops, i.e., closed curves in the Euclidean plane, with prescribed curvature $K$ (see \cite{BetCalGui, CalGui, GuiRol, KirLau, Mus, MusZud, NovVal}). Among them, we mention a couple of results more linked to our work: in \cite{KirLau} it is proved that if the prescribed curvature function $K$ is strictly monotone in a given direction, with non zero derivative, then no (immersed) $K$-loop exists. In \cite{Mus} it is shown that if the prescribed curvature function $K$ is positive, radially symmetric and is non-increasing as a function of the distance from the radial symmetry point, then any simple $K$-loop is in fact a circle.

In this work are interested in bounded, non round, possibly periodic motions with ``high energy'' for a class of mappings $B\colon\R^{2}\to\R$ which are radially symmetric (with respect to $u\in\R^{2}$) and with a behavior at infinity of the form
\begin{equation}
\label{eq:B-model}
B(u)\sim 1+\frac{A}{|u|^{\gamma}}\quad\text{as~~}|u|\to\infty~\!,\quad u\in\R^{2}~\!,
\end{equation}
with $\gamma>0$ and $A\in\R\setminus\{0\}$. Moreover, we consider a special class of curves, defined (in complex notation) as
$$
v_{\varepsilon,\rho}(t) = \rho e^{i\varepsilon t} + e^{it}
$$
with $\rho\gg 1$ and $0<|\varepsilon|\ll 1$. Such curves are obtained as a superposition of two motions: a counter-clockwise rotation with constant angular speed at a unit distance about a point which moves along a circumference of large radius $\rho$ with a slow angular speed $\varepsilon$ (see Fig.~1). The slow rotation is counter-clockwise if $\varepsilon>0$, or clockwise if $\varepsilon<0$. Observe that $\rho-1\le|v_{\varepsilon,\rho}(t)|\le\rho+1$ for every $t\in\R$, and that the (support of the) curve parametrized by $v_{\varepsilon,\rho}$ does not change after rotations of angles $\frac{2\pi\varepsilon}{1-\varepsilon}$ or $\frac{2\pi}{1-\varepsilon}$, because $v_{\varepsilon,\rho}$ satisfies
\begin{equation}
\label{eq:rotation}
e^{\frac{2\pi i\varepsilon}{1-\varepsilon}}v(t)=e^{\frac{2\pi i}{1-\varepsilon}}v(t)=v\left(t+\frac{2\pi}{1-\varepsilon}\right)~\!.
\end{equation}
In particular $v_{\varepsilon,\rho}$ is closed if and only if $\varepsilon$ is rational. 
\bigskip

\centerline{\begin{tabular}{ccc}
\begin{tikzpicture}[scale=0.7]
\begin{axis}[grid=major,
axis x line=middle,  axis y line=middle,
xtick=\empty, ytick=\empty,
xmin=-4.5,xmax=4.5,
ymin=-4.5,ymax=4.5,
axis equal image
]
\addplot [blue,domain=0:16*3.14,samples=120,smooth] 
({3*cos(deg((1/8)*x))+0.9* cos(deg(x))},{3*sin(deg((1/8)*x))+ 0.9*sin(deg(x))});
\end{axis}
\end{tikzpicture}
&
\begin{tikzpicture}[scale=0.7]
\begin{axis}[grid=major,
axis x line=middle,  axis y line=middle,
xtick=\empty, ytick=\empty,
xmin=-4.5,xmax=4.5,
ymin=-4.5,ymax=4.5,
axis equal image]
\addplot [blue,domain=0:26*3.14,samples=200,smooth] 
({3*cos(deg((2/13)*x))+0.9* cos(deg(x))},{3*sin(deg((2/13)*x))+ 0.9*sin(deg(x))});
\end{axis}
\end{tikzpicture}
&
\begin{tikzpicture}[scale=0.7]
\begin{axis}[grid=major,
axis x line=middle,  axis y line=middle,
xtick=\empty, ytick=\empty,
xmin=-4.5,xmax=4.5,
ymin=-4.5,ymax=4.5,
axis equal image]
\addplot [blue,domain=0:400*3.14,samples=1000,smooth] 
({3*cos(deg((exp(-2))*x))+0.9* cos(deg(x))},{3*sin(deg((exp(-2))*x))+ 0.9*sin(deg(x))});
\end{axis}
\end{tikzpicture}\\
\footnotesize{$\varepsilon=\frac{1}{8}$}
&
\footnotesize{$\varepsilon=\frac{2}{13}$}
&
\footnotesize{$\varepsilon\not\in\mathbb{Q}$}
\end{tabular}}
\bigskip

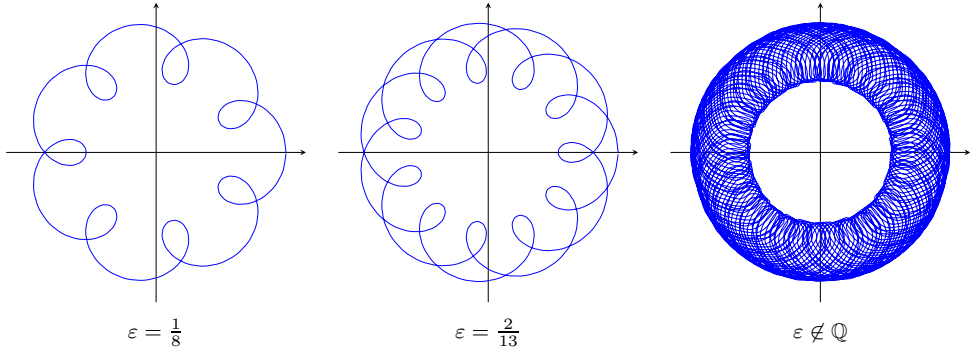
\centerline{\begin{minipage}[c]{0.82\linewidth}\footnotesize{Fig. 1 -- Curves parametrized by $v_{\varepsilon,\rho}$ for different values of $\varepsilon$. If $\varepsilon=\frac{m}{n}$ with $m,n\in\N$, $m<n$, then $v_{\varepsilon,\rho}$ is periodic of period $2n\pi$ and, for $\rho $ large enough, parametrizes a closed curve with $(n-m)$ curls and has a rotational symmetry of an angle $2\pi/(n-m)$. If $\varepsilon\not\in\mathbb{Q}$ then the range of $v_{\varepsilon,\rho}$ fills densely the annulus of radii $\rho-1$ and $\rho+1$.}
\end{minipage}}
\bigskip

Considering curves $v_{\varepsilon,\rho}$ is motivated by the observation that we look for solutions to (\ref{eq:curvature}) living far away, where $B$ is close to the constant value 1, in view of (\ref{eq:B-model}). Hence we expect to find solutions made by almost uniform circular motions of unit radius. Moreover, since we are interested in bounded trajectories and we deal with a planar problem characterized by radial symmetry, it seems reasonable to look for solutions of (\ref{eq:curvature}) as perturbations of $v_{\varepsilon,\rho}$ in the normal direction, i.e., of the form
$$
u=v_{\varepsilon,\rho}+\phi N_{\varepsilon,\rho}\quad\text{where~~}N_{\varepsilon,\rho}:=\frac{i\dot v_{\varepsilon,\rho}}{|\dot v_{\varepsilon,\rho}|}
$$
and $\phi$ is a scalar mapping. The discrete symmetry with respect to rotations of angles $\frac{2\pi\varepsilon}{1-\varepsilon}$ or $\frac{2\pi}{1-\varepsilon}$ is preserved asking $\phi$ to be periodic of period $\frac{2\pi}{1-\varepsilon}$. In fact, ve can prove the following result:

\begin{Theorem}\label{T:main}
Let $B \in C^1(\R^2; \R)$ be a radial function satisfying:
\begin{itemize}
\item[$(B_{1})$] 
$B(u)=1+A|u|^{-\gamma}+A_{1}|u|^{-\gamma_{1}}+o(|u|^{-\gamma-1})$ as $|u|\to\infty$, with $A,A_{1}\in\R$, $A\ne 0$, 
and $1<\gamma<\gamma_{1}$;
\item[$(B_{2})$]
$|\nabla B(u)|=O\left(|u|^{-\min\{\gamma+1,\gamma_{1}\}}\right)$ as $|u|\to\infty$.
\end{itemize}
If $A>0$ then there exist positive constants $\overline{\varepsilon}$, ${\mu}$, $a_{1}$, $a_{2}$ such that for every $\varepsilon\in(0,\overline{\varepsilon})$ there are $\rho_{\varepsilon}>0$ and a solution $u_{\varepsilon}$ of (\ref{eq:Lorentz3}) such that:
$$
a_{1}|\varepsilon|^{-\frac{1}{\gamma+2}}\le\rho_{\varepsilon}\le a_{2}|\varepsilon|^{-\frac{1}{\gamma+2}}\quad\text{and}\quad
\|u_{\varepsilon}-v_{\varepsilon,\rho_{\varepsilon}}\|_{C^{2}}\le{\mu}|\varepsilon|^{\frac{\gamma}{\gamma+2}}~\!.
$$
If $A<0$, then the same conclusion holds true with $\varepsilon\in(-\overline{\varepsilon},0)$. Moreover, in both cases, 
the curve parametrized by $u_{\varepsilon}$ is invariant with respect to rotations of angles $\frac{2\pi\varepsilon}{1-\varepsilon}$ and $\frac{2\pi}{1-\varepsilon}$ and is closed if $\varepsilon$ is rational.
\end{Theorem}

The assumption $\gamma>1$ as well as the conditions $(B_{1})$ and $(B_{2})$ precise the behavior of $B$ at infinity. As by-product of Theorem \ref{T:main}, we remark that in general the uniqueness result proved in \cite{Mus} and previously quoted is true just for simple loops. 

Going back to the problem of trajectories of a charged particle in a magnetic field with cylindrical symmetry, as a consequence of Theorem \ref{T:main}, we obtain:

\begin{Corollary}\label{C:main}
Let ${\bf{B}}\colon\R^{3}\to\R^{3}$ be a vector field of the form (\ref{eq:B3}) with $B\in C^1(\R^2; \R)$ radial function satisfying $(B_{1})$ and $(B_{2})$. Then there exists a family of motions for the Lorentz equation (\ref{eq:Lorentz}) whose projections on the horizontal plane, up to reparameterization in time, correspond to the solutions $u_{\varepsilon}$ given by Theorem \ref{T:main}. 
\end{Corollary}

The proof of Theorem \ref{T:main} is based on the Lyapunov-Schmidt reduction method, combined with a certain variational argument. Firstly, in order to get rid of the parameter $\varepsilon$, we define
\begin{equation}
\label{eq:u-epsilon-rho}
u_{\varepsilon,\rho}(t):=v_{\varepsilon,\rho}\left(\frac{t}{1-\varepsilon}\right)\quad\text{and}\quad
\n_{\varepsilon,\rho}(t):=N_{\varepsilon,\rho}\left(\frac{t}{1-\varepsilon}\right)
\end{equation}
and we look for solutions to (\ref{eq:curvature}) in the form
\begin{equation}
\label{eq:ansatz1}
u=u_{\varepsilon,\rho}+\varphi \n_{\varepsilon,\rho}
\end{equation}
with $\varphi$ periodic of period $2\pi$. Plugging (\ref{eq:ansatz1}) into (\ref{eq:curvature}) yields the following nonlinear problem for $\varphi$:
\begin{equation}
\label{eq:phi1}
\left\{\begin{array}{l}\varphi\colon\R\to\R~\!,~2\pi\text{-periodic}\\
\mathcal{K}_{\varepsilon,\rho}(\varphi)=B(u_{\varepsilon,\rho}+\varphi \n_{\varepsilon,\rho})\end{array}\right.
\end{equation}
where $\mathcal{K}_{\varepsilon,\rho}$ is the second order differential operator defined by
\begin{equation}
\label{eq:Ker}
\mathcal{K}_{\varepsilon,\rho}(\varphi):=\frac{i\frac{d}{dt}\left(u_{\varepsilon,\rho}+\varphi \n_{\varepsilon,\rho}\right)\cdot \frac{d^{2}}{dt^{2}}\left(u_{\varepsilon,\rho}+\varphi \n_{\varepsilon,\rho}\right)}{\left|\frac{d}{dt}\left(u_{\varepsilon,\rho}+\varphi \n_{\varepsilon,\rho}\right)\right|^{3}}~\!.
\end{equation}
Then, we identify the operator obtained by linearizing problem (\ref{eq:phi1}), in the limit $\varepsilon\to 0$. The resulting operator, acting among spaces of periodic mappings, turns out to have a two-dimensional kernel. After detaching its kernel, we can invert it and convert problem (\ref{eq:phi1}) into a fixed point problem, where the contraction principle can be applied. Actually, because of the two-dimensional kernel, we arrive to solve (\ref{eq:phi1}) for $|\varepsilon|$ small enough and for every $\rho$ sufficiently large in a suitable interval of values depending on $\varepsilon$, apart from a couple of Lagrange multipliers. This is performed in Section \ref{S:redux}. 

In order to remove the Lagrange multipliers, we exploit the variational nature of problem (\ref{eq:phi1}). In fact, (\ref{eq:phi1}) corresponds ultimately to the Euler-Lagrange equation of a suitable energy functional. One of the Lagrange multipliers can be eliminated by taking variations with respect to the parameter $\rho$. This leads to an equation for $\rho$ with respect to $\varepsilon$, which can be solved in correspondence of some $\rho_{\varepsilon}$ taking $\varepsilon>0$ if $A>0$ or $\varepsilon<0$ if $A<0$, as in the statement of Theorem \ref{T:main}. The other Lagrange multiplier comes from the rotational invariance and it disappears for free as soon as the first one vanishes. This part of the proof is accomplished in Section \ref{var}.

From the technical point of view, the argument used in the proof of Theorem \ref{T:main} shares some features with other works. In particular, we mention a higher dimensional version of the problem, about surfaces in $\R^3$ with prescribed mean curvature, shaped on collars of spheres arranged along a circumference with large radius, in order to construct a surface with genus 1. This kind of problem has been recently studied in \cite{CalMus, CalIacMus} for a class of mean curvature functions with a behavior at infinity like (\ref{eq:B-model}) but in the Euclidean 3-space.  We also quote the paper \cite{WeiYan}, where somehow similar techniques were already implemented to find multiple solutions to the nonlinear Schr\"oedinger equation with a radially symmetric potential again with a behavior like (\ref{eq:B-model}).

\section{The curvature operator}
\label{S:curvature-operator}

As described in the Introduction, the problem is redirected to the search of solutions to (\ref{eq:curvature}) in the form $u=u_{\varepsilon,\rho}+\varphi \n_{\varepsilon,\rho}$ where $u_{\varepsilon,\rho}(t)=\rho e^{i\varepsilon t/(1-\varepsilon)}+e^{it/(1-\varepsilon)}$, $\textbf{n}_{\varepsilon,\rho}(t)=i\dot u_{\varepsilon,\rho}(t)/|\dot u_{\varepsilon,\rho}(t)|$ and the unknown is the function $\varphi\colon\R\to\R$ which is asked to be $2\pi$-periodic. Hence the problem for $u$, in turn, is converted into the following nonlinear equation for $\varphi$: 
\begin{equation}
\label{eq:phi}
\mathcal{K}_{\varepsilon,\rho}(\varphi)=B(u_{\varepsilon,\rho}+\varphi \n_{\varepsilon,\rho})
\end{equation}
where $\mathcal{K}_{\varepsilon,\rho}$ is the second order differential operator defined by (\ref{eq:Ker}) and $B\colon\R^{2}\to\R$ is a radially symmetric $C^{1}$ function satisfying $(B_{1})$ and $(B_{2})$. Notice that if $u=u_{\varepsilon,\rho}+\varphi \n_{\varepsilon,\rho}$ and $\varphi$ is $2\pi$-periodic, then by \eqref{eq:rotation} and \eqref{eq:u-epsilon-rho}
\begin{equation}
\label{eq:rotation-u}
u(t+2\pi)=e^{\frac{2\pi i}{1-\varepsilon}}u(t)=e^{\frac{2\pi i\varepsilon}{1-\varepsilon}}u(t)~\!.
\end{equation}

In this Section we study the operator $\mathcal{K}_{\varepsilon,\rho}$, considered from the space of $2\pi$-periodic functions of class $C^{2}$, denoted $C^{2}(\R/2\pi\mathbb{Z})$, into the space of $2\pi$-periodic continuous functions, denoted $C^{0}(\R/2\pi\mathbb{Z})$. Such spaces are Banach spaces endowed with their natural norms. In fact $\mathcal{K}_{\varepsilon,\rho}$ turns out to be defined in 
\begin{equation}
\label{eq:N-domain}
\mathcal{N}:=\{\varphi\in C^{2}(\R/2\pi\mathbb{Z})~|~|\dot{u}_{\varepsilon,\rho}+\varphi\dot{\n}_{\varepsilon,\rho}+\dot{\varphi}\n_{\varepsilon,\rho}|>0\quad\text{on }\R\}~\!.
\end{equation}
\begin{Lemma}\label{L:curvature-operator}
The operator $\mathcal{K}_{\varepsilon,\rho}$ is of class $C^{\infty}$ from $\mathcal{N}$ into $C^{0}(\R/2\pi\mathbb{Z})$. In particular 
\begin{equation}
\label{eq:K'}
\mathcal{K}'_{\varepsilon,\rho}(\varphi)[\psi]=a_{\varepsilon,\rho}(\varphi)\psi''+b_{\varepsilon,\rho}(\varphi)\psi'+c_{\varepsilon,\rho}(\varphi)\psi\quad\forall\varphi\in\mathcal{N}~\!,~~\forall\psi\in C^{2}(\R/2\pi\mathbb{Z})~\!,
\end{equation}
where
\begin{equation}
\label{eq:coefficients-phi}
\begin{aligned}
&a_{\varepsilon, \rho }(\varphi) = \frac{i \dot u \cdot \n_{\varepsilon, \rho}}{|\dot u|^3},\\
&b_{\varepsilon, \rho}(\varphi) = \frac{2i \dot u \cdot \dot \n_{\varepsilon, \rho} - i \ddot u \cdot \n_{\varepsilon, \rho}}{|\dot u|^3}- \frac{3(i \dot u \cdot \ddot u)(\dot u \cdot \n_{\varepsilon, \rho})}{|\dot u |^5},\\
&c_{\varepsilon, \rho}(\varphi) = \frac{i \dot u \cdot \ddot \n_{\varepsilon, \rho}- i \ddot u \cdot \dot \n_{\varepsilon, \rho}}{|\dot u|^3} - \frac{3(i \dot u \cdot \ddot u)(\dot u \cdot \dot \n_{\varepsilon, \rho})}{|\dot u |^5},
\end{aligned}
\end{equation}
and $u={u}_{\varepsilon,\rho}+\varphi{\n}_{\varepsilon,\rho}$.
\end{Lemma}

\Proof
The operators $\varphi\mapsto i\frac{\de}{\de t}({u}_{\varepsilon,\rho}+\varphi{\n}_{\varepsilon,\rho})$ and $\varphi\mapsto\frac{\de^{2}}{\de t^{2}}({u}_{\varepsilon,\rho}+\varphi{\n}_{\varepsilon,\rho})$ are bounded linear operators from $C^{2}(\R/2\pi\mathbb{Z})$ into $C^{0}(\R/2\pi\mathbb{Z},\R^{3})$. Hence they are of class $C^{\infty}$. As the mapping $(v,w)\mapsto\frac{v\cdot w}{|v|^{3}}$ is of class $C^{\infty}$ in $(\R^{2}\setminus\{0\})\times\R^{2}$, the regularity of $\K_{\varepsilon,\rho}$ as well as the expression of $\K_{\varepsilon,\rho}'$ follows by standard differential calculus in Banach spaces. 
\QED

The differential operator $\Le_{0}\colon C^{2}(\R/2\pi\mathbb{Z})\to C^{0}(\R/2\pi\mathbb{Z})$ defined as
\begin{equation}
\label{linfop}
\Le_{0}\varphi := \ddot\varphi + \varphi
\end{equation}
will play a key role in the argument, because of the following fact. 

\begin{Lemma}\label{L:L0}
Fixing $0<a_{1}\le a_{2}$ and $\delta\in(0,1)$, one has that $\mathcal{K}_{\varepsilon,\rho}'(0)\to\Le_{0}$ as $\varepsilon\to 0$, uniformly with respect to $\rho\in S_{\varepsilon}:=[a_{1}|\varepsilon|^{-\delta}, a_{2}|\varepsilon|^{-\delta}]$, in the space of bounded linear operators from $C^{2}(\R/2\pi\mathbb{Z})$ in $C^{0}(\R/2\pi\mathbb{Z})$. More precisely, for $|\varepsilon|\ne 0$ small enough, and for every $\rho\in S_{\varepsilon}$ one has 
\begin{equation}\label{Lnest}
\|(\mathcal{K}'_{\varepsilon,\rho}(0)-\Le_{0})\varphi\|_{C^0}\le C|\varepsilon| \rho\|\varphi\|_{C^{2}}\quad\forall \varphi \in C^{2}(\R/2\pi\mathbb{Z})
\end{equation}
where $C$ depends only on $a_{1}$, $a_{2}$, and $\delta$.
\end{Lemma}

\Proof
Since
\begin{equation}\label{formulozze0}
\begin{split}
\dot \n_{\varepsilon,\rho} &= \frac{i \ddot u_{\varepsilon,\rho} }{|\dot u_{\varepsilon,\rho}|} -  \frac{( \dot u_{\varepsilon,\rho} \cdot \ddot u_{\varepsilon,\rho})  i \dot u_{\varepsilon,\rho}}{|\dot u_{\varepsilon,\rho}|^3}\\
\ddot \n_{\varepsilon,\rho} &= \frac{i \dddot u_{\varepsilon,\rho}}{|\dot u_{\varepsilon,\rho}|} -\frac{|\ddot u_{\varepsilon,\rho}|^2 i \dot u_{\varepsilon,\rho} + (\dot u_{\varepsilon,\rho} \cdot \dddot u_{\varepsilon,\rho})i \dot u_{\varepsilon,\rho} + 2(\dot u_{\varepsilon,\rho} \cdot \ddot u_{\varepsilon,\rho})i \ddot u_{\varepsilon,\rho}}{|\dot u_{\varepsilon,\rho}|^3} + 3 \frac{(\dot u_{\varepsilon,\rho} \cdot \ddot u_{\varepsilon,\rho})^2 i \dot u_{\varepsilon,\rho}}{|\dot u_{\varepsilon,\rho}|^5} 
\end{split}
\end{equation}
one has that
\begin{equation}\label{formulozze}
\begin{aligned}
&
i \dot u_{\varepsilon,\rho} \cdot \n_{\varepsilon,\rho} = |\dot u_{\varepsilon,\rho}| 
&~~~~~
&
\dot u_{\varepsilon,\rho} \cdot \n_{\varepsilon,\rho} = 0
&~~~~~
&
i \dot u_{\varepsilon,\rho} \cdot \ddot \n_{\varepsilon,\rho} = \frac{(\dot u_{\varepsilon,\rho} \cdot \ddot u_{\varepsilon,\rho})^2}{|\dot u_{\varepsilon,\rho}|^3} - \frac{|\ddot u_{\varepsilon,\rho}|^2}{|\dot u_{\varepsilon,\rho}|}
\\
&
\dot u_{\varepsilon,\rho} \cdot \dot \n_{\varepsilon,\rho} = -\frac{i \dot u_{\varepsilon,\rho} \cdot \ddot u_{\varepsilon,\rho}}{|\dot u_{\varepsilon,\rho}|}
&~~~~~
&
i \dot u_{\varepsilon,\rho} \cdot \dot \n_{\varepsilon,\rho} = 0 
&~~~~~
&
i\ddot u_{\varepsilon,\rho} \cdot \n_{\varepsilon,\rho} = \frac{\dot u_{\varepsilon,\rho} \cdot \ddot u_{\varepsilon,\rho}}{|\dot u_{\varepsilon,\rho}|}
\\
&
i \n_{\varepsilon,\rho}\cdot \dot \n_{\varepsilon,\rho} = \frac{i \dot u_{\varepsilon,\rho}\cdot \ddot u_{\varepsilon,\rho}}{|\dot u_{\varepsilon,\rho}|^2}
&~~~~~
&
\n_{\varepsilon,\rho}\cdot \dot \n_{\varepsilon,\rho} = 0
&~~~~~
&
i\ddot u_{\varepsilon,\rho} \cdot \dot \n_{\varepsilon,\rho}  =  \frac{|\ddot u_{\varepsilon,\rho}|^2}{|\dot u_{\varepsilon,\rho}|} - \frac{(\dot u_{\varepsilon,\rho} \cdot \ddot u_{\varepsilon,\rho})^2}{|\dot u_{\varepsilon,\rho}|^3}~\!.
\end{aligned}
\end{equation}
By (\ref{eq:K'}), (\ref{eq:coefficients-phi}) and (\ref{formulozze}) one infers that 
\begin{equation}
\label{eq:K'0}
\mathcal{K}_{\varepsilon,\rho}'(0)[\varphi]=a^0_{\varepsilon,\rho}\ddot\varphi+b^0_{\varepsilon,\rho}\dot\varphi+
c^0_{\varepsilon,\rho}\varphi
\end{equation}
where
$$
a^{0}_{\varepsilon,\rho} := \frac{1}{|\dot u_{\varepsilon,\rho}|^2}~\!,\quad b^{0}_{\varepsilon,\rho} := -\frac{\dot u_{\varepsilon,\rho}\cdot \ddot u_{\varepsilon,\rho}}{|\dot u_{\varepsilon,\rho}|^4}~\!,\quad c^{0}_{\varepsilon,\rho} := \frac{2(\dot u_{\varepsilon,\rho} \cdot \ddot u_{\varepsilon,\rho})^2 - 2|\ddot u_{\varepsilon,\rho}|^2 |\dot u_{\varepsilon,\rho}|^2  + 3(i \dot u_{\varepsilon,\rho} \cdot \ddot u_{\varepsilon,\rho})^2}{|\dot u_{\varepsilon,\rho}|^6}~\!.
$$
Fix $0<a_{1}\le a_{2}$ and $\delta\in(0,1)$.  Taking into account that 
$$
\begin{array}{ll}
|\dot u_{\varepsilon,\rho}|^{2}=(1-\varepsilon)^{-2}\left[1+2\varepsilon\rho\cos t+(\varepsilon\rho)^{2}\right]
~~~&
|\ddot u_{\varepsilon,\rho}|^{2}=(1-\varepsilon)^{-4}\left[1+2\varepsilon^{2}\rho\cos t+\varepsilon^{4}\rho^{2}\right]
\\
\dot u_{\varepsilon,\rho}\cdot\ddot u_{\varepsilon,\rho}=-(1-\varepsilon)^{-2}\varepsilon\rho\sin t
~~~&
i\dot u_{\varepsilon,\rho}\cdot\ddot u_{\varepsilon,\rho}=(1-\varepsilon)^{-3}\left[1+(\varepsilon^{2}\rho-\varepsilon\rho)\cos t-\varepsilon^{3}\rho^{2}\right],
\end{array}
$$
for $|\varepsilon|$ sufficiently small and non zero, $\rho \in S_{\varepsilon}$ and $\sigma\in\R$ the following estimates hold:  
\begin{equation}\label{stimeupic}
\begin{aligned}
&\left \||\dot u_{\varepsilon,\rho}|^\sigma- (1-\varepsilon)^{-\sigma}\right\|_{C^{0}} \le C|\varepsilon|\rho~\!, 
&~~~
&\left \||\ddot u_{\varepsilon,\rho}|^\sigma- (1-\varepsilon)^{-2\sigma}\right\|_{C^{0}} \leq C|\varepsilon|^2 \rho~\!, \\
&\left\|\dot u_{\varepsilon,\rho}\cdot \ddot u_{\varepsilon,\rho}\right\|_{C^{0}}\leq C|\varepsilon|\rho~\!,&~~~ &\left\|i \dot u_{\varepsilon,\rho}\cdot \ddot u_{\varepsilon,\rho}-(1-\varepsilon)^{-3}\right\|_{C^{0}} \leq C|\varepsilon|\rho~\!, \\
&\left \||\dddot u_{\varepsilon,\rho}|^\sigma- (1-\varepsilon)^{-3\sigma}\right\|_{C^{0}} \leq C|\varepsilon|^3 \rho~\!, &~~~ & \left\|\dot u_{\varepsilon,\rho}\cdot \dddot u_{\varepsilon_\rho}-(1-\varepsilon)^{-4}\right\|_{C^{0}}\le C|\varepsilon|\rho~\!,
\end{aligned}
\end{equation}
where $C$ depends only on $a_1, a_2, \delta$. By means of \eqref{stimeupic} one obtains that
$$
\|a^{0}_{\varepsilon,\rho} - 1\|_{C^{0}} \le C|\varepsilon|\rho~\!, \quad
\|b^{0}_{\varepsilon,\rho}\|_{C^{0}}\le C|\varepsilon|\rho~\!, \quad
\|c^{0}_{\varepsilon,\rho}-1\|_{C^{0}}\le C|\varepsilon|\rho~\!,
$$
The conclusion follows from these estimates and from (\ref{eq:K'0}).
\QED

Thanks to Lemma \ref{L:L0}, the operator $\mathcal{L}_{0}$ defined by \eqref{linfop} will constitute the leading term of the curvature operator, in the limit $\varepsilon\to 0$. It is well known that
$$
\mathrm{ker}(\mathcal{L}_{0})=\{\varphi\in C^{2}(\R/2\pi\mathbb{Z})~|~\Le_{0}\varphi=0\}=\mathrm{span}\{\cos t,~\!\sin t\}~\!.
$$
Moreover, for an arbitrary continuous bounded function $f\colon\R\to\R$, the equation $\Le_{0} \varphi = f$ admits a two-parameter family of global solutions, which turn out to be periodic if and only if
\begin{equation}\label{perpendicul}
\int_0^{2\pi}f(t) \cos t \de t = \int_0^{2\pi}f(t) \sin t \de t = 0~\!.
\end{equation}
Justified by that, we introduce the decomposition
\[
\begin{aligned}
&C^{2}(\R/2\pi\mathbb{Z}) = \text{span}\{\cos t, \sin t\} + X^{\perp}\quad\text{where~~}X^{\perp}:=\big\{f\in C^{2}(\R/2\pi\mathbb{Z}) ~|~\text{(\ref{perpendicul}) holds}\big\}\!~, \\
&C^{0}(\R/2\pi\mathbb{Z}) = \text{span}\{\cos t, \sin t\} + Y^{\perp}
\quad\text{where~~}Y^{\perp}:=\big\{f\in C^{0}(\R/2\pi\mathbb{Z}) ~|~\text{(\ref{perpendicul}) holds}\big\}\!~.
\end{aligned}
\]
Hence the operator $\mathcal{L}_{0}$ is a bijection from $X^\perp$ onto $Y^\perp$, it is linear and continuous and, by well known facts, also $\mathcal{L}_{0}^{-1}\colon Y^{\perp}\to X^{\perp}$ is so. Therefore there exists a constant $C_{0}>0$ such that 
\begin{equation}\label{linearestimate}
\|\varphi\|_{C^{2}} \leq C_{0}\|f\|_{C^{0}} \quad \forall f \in Y^\perp
\end{equation}
where $\varphi$ is the unique solution of $\Le_{0} \varphi = f$ in $X^{\perp}$. 
\medskip

In the remaining part of the work we aim to prove the following result.
\begin{Theorem}\label{mainteo}
Let $B \in C^1(\R^2; \R)$ be a radial function satisfying $(B_{1})$ and $(B_{2})$. There exist positive constants $\overline\varepsilon$, $a_{1}$, $a_{2}$, $\mu_{0}$ such that if $A > 0$, then for any $\varepsilon\in(0,\overline\varepsilon)$ there exist $\rho_\varepsilon\in\big[a_{1}|\varepsilon|^{-\frac{1}{\gamma+2}},a_{2}|\varepsilon|^{-\frac{1}{\gamma+2}}\big]$ and $\varphi_\varepsilon \in C^{2}(\R/2\pi\mathbb{Z})$ solving (\ref{eq:phi}) and such that $\|\varphi_\varepsilon\|_{C^{2}}\le\mu_0 |\varepsilon|^{\frac{\gamma}{\gamma+2}}$. If $A<0$ the same conclusion holds true but for $\varepsilon\in(-\overline\varepsilon,0)$. 
\end{Theorem}

Noting that, by \eqref{formulozze0} and \eqref{stimeupic}, $\sup_{0<|\varepsilon|<\overline\varepsilon}\|\n_{\varepsilon,\rho_{\varepsilon}}\|_{C^{2}}<\infty$, Theorem \ref{T:main} is an immediate consequence of Theorem \ref{mainteo}, whereas Corollary \ref{C:main} easily follows, by the next elementary result:

\begin{Lemma}\label{L:curvature}
If $u\colon\R\to\R^{2}$ is a classical solution of (\ref{eq:curvature}), 
then there exists an increasing diffeomorphism $g$ of $\R$ onto $\R$ such that $u\circ g$ solves (\ref{eq:Lorentz3}).
\end{Lemma}
\Proof
Let $u(t)$ be a solution to (\ref{eq:curvature}). The arc length $\ell(t)=\int_{0}^{t}|\dot u(\tau)|~\!d\tau$ is a $C^{1}$, invertible mapping and denoting $g(s)$ its inverse, the function $v=u\circ g$ satisfies 
$$
v'(s)=\frac{\dot{u}(g(s))}{|\dot{u}(g(s))|}~\!,\quad v''(s)=\frac{\ddot{u}(g(s))}{|\dot{u}(g(s))|^{2}}-\frac{\dot{u}(g(s))\cdot\ddot{u}(g(s))}{|\dot{u}(g(s))|^{4}}\dot{u}(g(s))~\!.
$$
In particular $|v'|\equiv 1$ and $v''\cdot v'=0$. Then $v''$ is parallel to $iv'$ and, in view of (\ref{eq:curvature}), $v$ solves (\ref{eq:Lorentz3}).
\QED

\section{The finite-dimensional reduction}\label{S:redux}

In this Section we tackle equation (\ref{eq:phi}) and we prove that it admits solutions for $|\varepsilon|$ small enough and non zero and for every $\rho$ sufficiently large, up to an extra term which is a linear combination of $\sin t$ and $\cos t$. More precisely, the following result holds true:
 
\begin{Theorem}\label{T:reduced}
Fix $0<a_1\leq a_2$ and $\delta \in (0,1)$. Then there exists $\varepsilon_{0}>0$ depending on $a_1$, $a_2$ and $\delta$ such that for every $0<|\varepsilon|\le \varepsilon_{0}$ and for every $\rho\in\big[a_1|\varepsilon|^{-\delta},a_2|\varepsilon|^{-\delta}\big]$ there exists $\varphi_{\varepsilon,\rho}\in C^{2}(\R/2\pi\mathbb{Z})$ solving 
\begin{equation}\label{eqconmolt}
\mathcal{K}_{\varepsilon,\rho}(\varphi_{\varepsilon,\rho})-B(u_{\varepsilon,\rho}+\varphi_{\varepsilon,\rho} \n_{\varepsilon,\rho})=\lambda^{1}_{\varepsilon,\rho}\cos t+\lambda^{2}_{\varepsilon,\rho}\sin t
\end{equation}
where 
\begin{equation}{\label{lambrappr}}
\begin{split}
\lambda^{1}_{\varepsilon,\rho}&=\frac{1}{\pi}\int_{0}^{2\pi}\left[\mathcal{K}_{\varepsilon,\rho}(\varphi_{\varepsilon,\rho})-B(u_{\varepsilon,\rho}+\varphi_{\varepsilon,\rho} n_{\varepsilon,\rho})\right]\cos t~\!dt~\!,\\
\lambda^{2}_{\varepsilon,\rho}&=\frac{1}{\pi}\int_{0}^{2\pi}\left[\mathcal{K}_{\varepsilon,\rho}(\varphi_{\varepsilon,\rho})-B(u_{\varepsilon,\rho}+\varphi_{\varepsilon,\rho} n_{\varepsilon,\rho})\right]\sin t~\!dt~\!.
\end{split}
\end{equation}
Moreover
\begin{equation}\label{estimphi}
\|\varphi_{\varepsilon,\rho}\|_{C^{2}}\le\mu_{0}|\varepsilon|^{\widehat\gamma}
\end{equation}
where $\widehat\gamma=\min\{\gamma\delta,1-\delta\}$ and $\mu_{0}$ is a constant depending on $a_1$, $a_2$ and $\delta$ but not on $\varepsilon$ neither on $\rho$. 
\end{Theorem}
 
Let $0 < a_1\leq a_2$ and $\delta \in (0,1)$ be fixed and let us introduce some notation. For every $\varepsilon\ne 0$ we define
$$
S_\varepsilon:= \big[a_1|\varepsilon|^{-\delta}, a_2|\varepsilon|^{-\delta}\big]~\!.
$$
We will always take $\rho \in S_\varepsilon$. Then, we fix $\varepsilon_{1}\in\big(0,\frac{1}{2}\big]$ satisfying
\begin{equation}\label{nsotto}
2<a_1\varepsilon_{1}^{-\delta} \text{  and  } a_2\varepsilon_{1}^{1-\delta} <1~\!.
\end{equation}
Hence for every $\varepsilon\in[-\varepsilon_{1},\varepsilon_{1}]\setminus\{0\}$ one has that $S_\varepsilon \subset (2, |\varepsilon|^{-1})$. In particular we recover the following inequality, extensively used in the sequel: for $\varepsilon\in[-\varepsilon_{1},\varepsilon_{1}]\setminus\{0\}$ and $\rho \in S_\varepsilon$, it holds
\begin{equation}\label{passag17}
1 < C_1 |\varepsilon|^{-\delta} \leq |u_{\varepsilon,\rho}(t)| \leq C_2|\varepsilon|^{-\delta} \quad \forall t \in \R~\!,
\end{equation}
for some positive constants $C_1$ and $C_2$ depending only on $a_1$, $a_2$ and $\delta$. 

Let us introduce the operator $\mathcal{B}_{\varepsilon,\rho}\colon C^{2}(\R/2\pi\mathbb{Z})\to C^{0}(\R/2\pi\mathbb{Z})$ defined by
$$
\mathcal{B}_{\varepsilon,\rho}(\varphi):=B(u_{\varepsilon,\rho}+\varphi\n_{\varepsilon,\rho})~\!.
$$
The operator $\mathcal{B}_{\varepsilon,\rho}$ is well defined because $B$ is radial and \eqref{eq:rotation-u} holds. Moreover it inherits the same regularity of $B$, i.e., it is of class $C^{1}$. In particular
\begin{equation}\label{mathcal-B'}
\mathcal{B}'_{\varepsilon,\rho}(\varphi)[\psi]=\psi\nabla B(u_{\varepsilon,\rho}+\varphi\n_{\varepsilon,\rho})\cdot\n_{\varepsilon,\rho}\quad\forall\varphi,\psi\in C^{2}(\R/2\pi\mathbb{Z})~\!.
\end{equation}
Finally, we define the operator $\mathcal{F}_{\varepsilon,\rho}\colon\mathcal{N}\to C^{0}(\R/2\pi\mathbb{Z})$ by setting
\begin{equation}\label{funfF}
\mathcal{F}_{\varepsilon,\rho}(\varphi) :=\mathcal{B}_{\varepsilon,\rho}(\varphi)-\mathcal{K}_{\varepsilon,\rho}(\varphi)+\mathcal{L}_{0}\varphi~\!.
\end{equation}
In view of the regularity of $\mathcal{B}_{\varepsilon,\rho}$ and by Lemma \ref{L:curvature-operator}, $\mathcal{F}_{\varepsilon,\rho}$ is of class $C^{1}$ in $\mathcal{N}$. Moreover 
\begin{gather*}
\mathcal{F}_{\varepsilon,\rho}(0)=B(u_{\varepsilon,\rho})-\mathcal{K}(u_{\varepsilon,\rho})~\!,\\
\mathcal{F}'_{\varepsilon,\rho}(0)[\varphi]=\varphi\nabla B(u_{\varepsilon,\rho})\cdot\n_{\varepsilon,\rho}-\mathcal{K}'_{\varepsilon,\rho}(0)[\varphi]+\mathcal{L}_{0}\varphi~~\forall\varphi\in  C^{2}(\R/2\pi\mathbb{Z})
\end{gather*}
with $\mathcal{K}$ defined in \eqref{mathcal-K}. 

The following estimates will be useful. From now on, with the symbol $\|\cdot\|$ without subscript we will denote the norm in the space of bounded linear operators from $C^2(\R/2\pi\Z)$ into $C^0(\R/2\pi\Z)$.

\begin{Lemma}\label{L:reduced}
There exist $\varepsilon_{0}\in(0,\varepsilon_{1}]$ and $M_0>0$ such that for every $0<|\varepsilon|\le\varepsilon_{0}$ and $\rho \in S_\varepsilon$ it holds
\begin{gather}\label{error}
\|\mathcal{F}_{\varepsilon,\rho}(0)\|_{C^0}\le M_{0}|\varepsilon|^{\widehat\gamma}\\
\label{error-1}
\|\mathcal{F}'_{\varepsilon,\rho}(\varphi)\|\le\frac{1}{8C_{0}}\text{~~if~~}\|\varphi\|_{C^{2}}\le\mu_{0}|\varepsilon|^{\widehat\gamma}
\end{gather}
where $C_{0}$ is the constant in \eqref{linearestimate}, $\widehat\gamma:=\min\{\gamma\delta, 1-\delta\}$ and $\mu_{0}:=8C_{0}M_{0}$.
\end{Lemma}

\Proof
We have that
\begin{equation}
\label{error-bis}
\|\mathcal{F}_{\varepsilon,\rho}(0)\|_{C^0}\le
\|B(u_{\varepsilon,\rho})-1\|_{C^0}+\left\|\mathcal{K}(u_{\varepsilon,\rho})-1\right\|_{C^0}.
\end{equation}
Since, by $(B_{1})$, $B(v)=1+A|v|^{-\gamma}+o(|v|^{-\gamma})$ as $|v|\to\infty$, using also \eqref{passag17}, one can find $\varepsilon_{0}\in(0,\varepsilon_{1}]$ such that
\begin{equation}
\label{error-B}
\|B(u_{\varepsilon,\rho})-1\|_{C^0}\le 2|A|C_{1}^{-\gamma}|\varepsilon|^{\gamma\delta}\quad\forall\rho\in S_{\varepsilon}~\!,~~\forall\varepsilon\in[-\varepsilon_{0},\varepsilon_{0}]\setminus\{0\}~\!.
\end{equation}
On the other hand, by \eqref{stimeupic} we obtain
\begin{equation}\label{passag1948}
\left\|\mathcal{K}(u_{\varepsilon,\rho})-1\right\|_{C^0}\leq C|\varepsilon|^{1-\delta}\quad\forall\rho\in S_{\varepsilon}~\!,~~\forall\varepsilon\in[-\varepsilon_{0},\varepsilon_{0}]\setminus\{0\}
\end{equation}
where $C$ is a constant depending only on $a_{1},a_{2}$ and $\delta$. Hence \eqref{error-bis}--\eqref{passag1948} imply \eqref{error}, for some $M_{0}$ depending on $|A|$, $C_{1}$, and $C$. Let us show \eqref{error-1}. We have that 
\begin{equation}\label{F'}
\|\mathcal{F}'_{\varepsilon,\rho}(0)\|\le
\|\nabla B(u_{\varepsilon,\rho})\|_{C^{0}}+
\|\mathcal{K}'_{\varepsilon,\rho}(0)-\mathcal{L}_{0}\|
\end{equation}
By ($B_{2}$) and \eqref{stimeupic} we have that
\begin{equation}\label{F'1}
\|\nabla B(u_{\varepsilon,\rho})\|_{C^{0}}\le C|\varepsilon|^{\delta\min\{\gamma+1,\gamma_{1}\}}\quad\forall\rho\in S_{\varepsilon}~\!,~~\forall\varepsilon\in[-\varepsilon_{1},\varepsilon_{1}]\setminus\{0\}~\!.
\end{equation}
Moreover Lemma \ref{L:L0} yields that
\begin{equation}\label{F'2}
\|\mathcal{K}'_{\varepsilon,\rho}(0)-\mathcal{L}_{0}\|\le C|\varepsilon|\rho\le C|\varepsilon|^{1-\delta}\quad\forall\rho\in S_{\varepsilon}~\!,~~\forall\varepsilon\in[-\varepsilon_{1},\varepsilon_{1}]\setminus\{0\}~\!.
\end{equation}
By \eqref{F'}--\eqref{F'2} and taking a smaller $\varepsilon_{0}\in(0,\varepsilon_{1}]$, we obtain that
\begin{equation}
\label{F0-zero}
\sup_{\rho\in S_{\varepsilon}}\|\mathcal{F}'_{\varepsilon,\rho}(0)\|\to 0\quad\text{as~~}\varepsilon\to 0~\!.
\end{equation}
Now we show that
\begin{equation}
\label{Fclaim}
\sup_{\scriptstyle{\rho\in S_{\varepsilon}}\atop\scriptstyle{\varphi\in\mathcal{M}_{\varepsilon}}}\|\mathcal{F}'_{\varepsilon,\rho}(\varphi)-\mathcal{F}'_{\varepsilon,\rho}(0)\|\to 0\quad\text{as~~}\varepsilon\to 0~\!,
\end{equation}
where $\mathcal{M}_{\varepsilon}:=\{\varphi\in C^{2}(\R/2\pi\mathbb{Z})~|~\|\varphi\|_{C^{2}}\le \mu_{0}|\varepsilon|^{\widehat\gamma}\}$. 
Indeed, from \eqref{passag17}--\eqref{funfF} it follows that
\begin{equation}
\label{Fclaim1}
\|\mathcal{F}'_{\varepsilon,\rho}(\varphi)-\mathcal{F}'_{\varepsilon,\rho}(0)\|\le\|\nabla B(u_{\varepsilon,\rho}+\varphi\n_{\varepsilon,\rho})-\nabla B(u_{\varepsilon,\rho})\|_{C^{0}}+\|\mathcal{K}'_{\varepsilon,\rho}(\varphi)-\mathcal{K}'_{\varepsilon,\rho}(0)\|
\end{equation}
Moreover, from $(B_{2})$
$$
|\nabla B(u_{\varepsilon,\rho}+\varphi\n_{\varepsilon,\rho})|\le C|u_{\varepsilon,\rho}+\varphi\n_{\varepsilon,\rho})|^{-\min\{\gamma+1,\gamma_{1}\}}\le C(|u_{\varepsilon,\rho}|-|\varphi|)^{-\min\{\gamma+1,\gamma_{1}\}}
$$
and using \eqref{passag17} and the bound for $\varphi\in\mathcal{M}_{\varepsilon}$, we deduce that
$$
\|\nabla B(u_{\varepsilon,\rho}+\varphi\n_{\varepsilon,\rho})\|_{C^{0}}\le \frac{C}{(C_{1}|\varepsilon|^{-\delta}-\mu_{0}|\varepsilon|^{\widehat\gamma})^{\min\{\gamma+1,\gamma_{1}\}}}~\!.
$$
This estimate together with \eqref{F'1} yields
\begin{equation}
\label{Bclaim}
\sup_{\scriptstyle{\rho\in S_{\varepsilon}}\atop\scriptstyle{\varphi\in\mathcal{M}_{\varepsilon}}}\|\nabla B(u_{\varepsilon,\rho}+\varphi\n_{\varepsilon,\rho})-B(u_{\varepsilon,\rho})\|_{C^{0}}\to 0\quad\text{as~~}\varepsilon\to 0~\!
\end{equation}
Now we show that
\begin{equation}
\label{abc0}
\sup_{\scriptstyle{\rho\in S_{\varepsilon}}\atop\scriptstyle{\varphi\in\mathcal{M}_{\varepsilon}}}\|\mathcal{K}'_{\varepsilon,\rho}(\varphi)-\mathcal{K}'_{\varepsilon,\rho}(0)\|\le C|\varepsilon|^{\widehat\gamma}\quad\forall\varepsilon\in[-\varepsilon_{1},\varepsilon_{1}]\setminus\{0\}~\!.
\end{equation}
By \eqref{eq:K'}, ve have that
\begin{equation}
\label{abc1}
\|\mathcal{K}'_{\varepsilon,\rho}(\varphi)-\mathcal{K}'_{\varepsilon,\rho}(0)\|\le\|a_{\varepsilon,\rho}(\varphi)-a_{\varepsilon,\rho}(0)\|_{C^{0}}+\|b_{\varepsilon,\rho}(\varphi)-b_{\varepsilon,\rho}(0)\|_{C^{0}}+\|c_{\varepsilon,\rho}(\varphi)-c_{\varepsilon,\rho}(0)\|_{C^{0}}
\end{equation}
where $a_{\varepsilon,\rho}(\varphi)$, $b_{\varepsilon,\rho}(\varphi)$ and $c_{\varepsilon,\rho}(\varphi)$ are defined in \eqref{eq:coefficients-phi}. Setting $\Phi_{\varepsilon,\rho}=\varphi\n_{\varepsilon,\rho}$, we write
\begin{equation*}
\begin{aligned}
a_{\varepsilon,\rho}(\varphi)-a_{\varepsilon,\rho}(0)&=i\big[\widetilde{a}(\dot u_{\varepsilon,\rho}+\dot\Phi_{\varepsilon,\rho})-\widetilde{a}(\dot u_{\varepsilon,\rho})\big]\cdot\n_{\varepsilon,\rho},\\
b_{\varepsilon,\rho}(\varphi)-b_{\varepsilon,\rho}(0)&=2i\big[\widetilde{a}(\dot u_{\varepsilon,\rho}+\dot\Phi_{\varepsilon,\rho})-\widetilde{a}(\dot u_{\varepsilon,\rho})\big]\!\cdot\!\dot\n_{\varepsilon,\rho}-i\big[\widetilde{b}(\dot u_{\varepsilon,\rho}+\dot\Phi_{\varepsilon,\rho},\ddot u_{\varepsilon,\rho}+\ddot\Phi_{\varepsilon,\rho})-\widetilde{b}(\dot u_{\varepsilon,\rho},\ddot u_{\varepsilon,\rho})\big]\!\cdot\!\n_{\varepsilon,\rho}\\
&\quad-3 \big[\widetilde{c}(\dot u_{\varepsilon,\rho}+\dot\Phi_{\varepsilon,\rho},\ddot u_{\varepsilon,\rho}+\ddot\Phi_{\varepsilon,\rho})-\widetilde{c}(\dot u_{\varepsilon, \rho}, \ddot u_{\varepsilon,\rho})\big]\!\cdot\!\n_{\varepsilon,\rho},\\
c_{\varepsilon,\rho}(\varphi)-c_{\varepsilon,\rho}(0)&=i\big[\widetilde{a}(\dot u_{\varepsilon,\rho}+\dot\Phi_{\varepsilon,\rho})-\widetilde{a}(\dot u_{\varepsilon,\rho})\big]\!\cdot\!\ddot\n_{\varepsilon,\rho}-i\big[\widetilde{b}(\dot u_{\varepsilon,\rho}+\dot\Phi_{\varepsilon,\rho},\ddot u_{\varepsilon,\rho}+\ddot\Phi_{\varepsilon,\rho})-\widetilde{b}(\dot u_{\varepsilon,\rho},\ddot u_{\varepsilon,\rho})\big]\!\cdot\! \dot \n_{\varepsilon,\rho}\\
&\quad-3 \big[\widetilde{c}(\dot u_{\varepsilon,\rho}+\dot\Phi_{\varepsilon,\rho},\ddot u_{\varepsilon,\rho}+\ddot\Phi_{\varepsilon,\rho})-\widetilde{c}(\dot u_{\varepsilon, \rho}, \ddot u_{\varepsilon,\rho})\big]\!\cdot\!\dot\n_{\varepsilon,\rho},\\
\end{aligned}
\end{equation*}
where
$$
\widetilde{a}(v)=\frac{v}{|v|^{3}}~\!,~~
\widetilde{b}(v,w)=\frac{w}{|v|^{3}}~\!,~~
\widetilde{c}(v,w)=\frac{(iv\cdot w)v}{|v|^{5}}\quad\forall(v,w)\in(\R^{2}\setminus\{0\})\times\R^{2}~\!.
$$
By \eqref{formulozze0} and \eqref{stimeupic}, there exists a constant $C>0$ independent of $\varepsilon$ and $\rho$, such that $\|\n_{\varepsilon,\rho}\|_{C^{2}}\le C$ for every $\varepsilon\in[-\varepsilon_{0},\varepsilon_{0}]\setminus\{0\}$ and for every $\rho\in S_{\varepsilon}$. From this it also follows that 
$$
\|\varphi\n_{\varepsilon,\rho}\|_{C^{2}}\le C|\varepsilon|^{\widehat\gamma}\quad\forall \varepsilon\in[-\varepsilon_{0},\varepsilon_{0}]\setminus\{0\}~\!,~~\forall\rho\in S_{\varepsilon}~\!,~~\forall\varphi\in\mathcal{M}_{\varepsilon}~\!.
$$
Moreover, by \eqref{stimeupic}, there exists a compact set $K\subset(\R^{2}\setminus\{0\})\times\R^{2}$ such that 
$(\dot u_{\varepsilon,\rho}(t)+s\dot\Phi_{\varepsilon,\rho}(t),\ddot u_{\varepsilon,\rho}(t)+s\ddot\Phi_{\varepsilon,\rho}(t))\in K$ for all $t\in\R$, $s\in[0,1]$, $\varepsilon\in[-\varepsilon_{0},\varepsilon_{0}]\setminus\{0\}$, $\rho\in S_{\varepsilon}$ and $\varphi\in\mathcal{M}_{\varepsilon}$. Using these facts as well as the regularity of the functions $\widetilde{a}$, $\widetilde{b}$ and $\widetilde{c}$, by the mean value theorem, we infer that
\begin{equation}
\label{abc2}
\|a_{\varepsilon,\rho}(\varphi)-a_{\varepsilon,\rho}(0)\|_{C^{0}}\le C|\varepsilon|^{\widehat\gamma}~\!,~~~
\|b_{\varepsilon,\rho}(\varphi)-b_{\varepsilon,\rho}(0)\|_{C^{0}}\le C|\varepsilon|^{\widehat\gamma}~\!,~~~
\|c_{\varepsilon,\rho}(\varphi)-c_{\varepsilon,\rho}(0)\|_{C^{0}}\le C|\varepsilon|^{\widehat\gamma}
\end{equation}
for every $\varepsilon\in[-\varepsilon_{0},\varepsilon_{0}]\setminus\{0\}$, $\rho\in S_{\varepsilon}$ and $\varphi\in\mathcal{M}_{\varepsilon}$, and $C$ positive constant independent of $\varepsilon$, $\rho$ and $\varphi$. Hence \eqref{abc1} and \eqref{abc2} imply \eqref{abc0}, and \eqref{Fclaim} follows from \eqref{Fclaim1}--\eqref{abc0}. Finally \eqref{error-1} is a consequence of \eqref{F0-zero} and \eqref{Fclaim}, for a possibly smaller $\varepsilon_{0}$.
\QED

\noindent
\emph{Proof of Theorem \ref{T:reduced}.}
By \eqref{funfF}, \eqref{eq:phi} is equivalent to 
\begin{equation}\label{linearform}
\mathcal{L}_{0}\varphi=\mathcal{F}_{\varepsilon,\rho}(\varphi)~\!. 
\end{equation}
We aim to rewrite \eqref{linearform} as a fixed point problem in $C^{2}(\R/2\pi\mathbb{Z})$. Since $\textrm{range}(\Le_{0})= Y^\perp$, we consider the projection of $\F_{\varepsilon,\rho}(\varphi)$ on $Y^\perp$, given by 
\[
\hat \F_{\varepsilon,\rho}(\varphi) := \F_{\varepsilon,\rho}(\varphi) +\lambda^1_{\varepsilon,\rho}(\varphi) \cos t + \lambda^2_{\varepsilon,\rho}(\varphi)  \sin t, 
\]
where
$$
\lambda^1_{\varepsilon,\rho}(\varphi)  := - \frac{1}{\pi}\int_0^{2\pi}\F_{\varepsilon,\rho}(\varphi)\cos t \de t~\!, \quad \lambda^2_{\varepsilon,\rho}(\varphi)  := - \frac{1}{\pi}\int_0^{2\pi}\F_{\varepsilon,\rho}(\varphi)\sin t \de t~\!.
$$
Observe that
\begin{equation}
\label{Fhat}
\begin{array}{c}
\|\hat \F_{\varepsilon,\rho}(\varphi)\|_{C^{0}}\le 4\|\F_{\varepsilon,\rho}(\varphi)\|_{C^{0}}\quad\forall\varphi\in\mathcal{N}\vspace{4pt}\\ 
\|\hat \F_{\varepsilon,\rho}(\varphi_{1})-\hat \F_{\varepsilon,\rho}(\varphi_{2})\|_{C^{0}}\le 4\|\F_{\varepsilon,\rho}(\varphi_{1})-\F_{\varepsilon,\rho}(\varphi_{2})\|_{C^{0}}\quad\forall\varphi_{1},\varphi_{2}\in\mathcal{N}~\!.
\end{array}
\end{equation}
Then we set
$$
\Q_{\varepsilon,\rho}:= \Le_{0}^{-1}\circ \hat \F_{\varepsilon,\rho}~\!,
$$
so that if $\varphi\in C^{2}(\R/2\pi\mathbb{Z})$ solves 
\begin{equation}\label{fixpopoint}
\varphi = \Q_{\varepsilon,\rho}(\varphi)
\end{equation}
then it satisfies $\mathcal{L}_{0}\varphi={\hat \F}_{\varepsilon,\rho}(\varphi)$ and thus \eqref{eqconmolt}. 
We can solve (\ref{fixpopoint}) in a suitable neighborhood $\mathcal{N}_{\varepsilon}$ of $0$ in $Y^{\perp}$, whose size is determined by the estimate of the ``error'' $\mathcal{F}_{\varepsilon}(0)$ according to Lemma \ref{L:reduced}. More precisely, we set
$$
\mathcal{N}_{\varepsilon}:=\{\varphi\in X^{\perp}\cap\mathcal{N}~|~\|\varphi\|_{C^{2}}\le\mu_{0}|\varepsilon|^{\widehat\gamma}\}
$$
where $\mathcal{N}$ is defined in \eqref{eq:N-domain} and $\widehat\gamma$ and $\mu_{0}$ are given by Lemma \ref{L:reduced}. By construction, $\hat \F_{\varepsilon,\rho}(\varphi) \in Y^\perp$. Moreover $\mathcal{L}_{0}^{-1}\colon Y^{\bot}\to X^{\bot}$. Therefore, $\Q_{\varepsilon,\rho}(\varphi)\in X^\perp$ for every $\varphi\in\mathcal{N}_{\varepsilon}$ and, thanks to \eqref{linearestimate} and \eqref{Fhat},
\begin{gather}
\label{M1}
\|\Q_{\varepsilon,\rho}(\varphi)\|_{C^2} \leq 4C_{0}\|\F_{\varepsilon,\rho}(\varphi)\|_{C^0}\quad \forall \varphi \in \mathcal{N}_{\varepsilon}\\
\label{Lip}
\|\Q_{\varepsilon,\rho}(\varphi_2)- \Q_{\varepsilon,\rho}(\varphi_1)\|_{C^2}\leq 4C_{0}\|\F_{\varepsilon,\rho}(\varphi_1) - \F_{\varepsilon,\rho}(\varphi_2)\|_{C^{0}}\quad \forall \varphi_1, \varphi_2 \in \mathcal{N}_{\varepsilon}
\end{gather}
where $C_{0}$ is the constant in (\ref{linearestimate}). If $\varphi_1, \varphi_2 \in \mathcal{N}_{\varepsilon}$ then, also $s\varphi_1+(1-s)\varphi_2\in\mathcal{N}_{\varepsilon}$ for every $s\in[0,1]$ and 
\begin{equation}
\label{FLip}
\|\F_{\varepsilon,\rho}(\varphi_1) - \F_{\varepsilon,\rho}(\varphi_2)\|_{C^{0}}\le\max_{s\in[0,1]}\|\F'_{\varepsilon,\rho}(s\varphi_1+(1-s)\varphi_2)\|\|\varphi_1-\varphi_2\|_{C^{2}}\le\frac{1}{8C_{0}}\|\varphi_1-\varphi_2\|_{C^{2}}
\end{equation}
thanks to \eqref{error-1}. Hence, \eqref{Lip} and \eqref{FLip} imply that $\Q_{\varepsilon,\rho}$ is a contraction in $\mathcal{N}_{\varepsilon}$. Moreover, if $\varphi\in \mathcal{N}_{\varepsilon}$ then, by \eqref{error} and \eqref{FLip}, and recalling that $\mu_{0}=8C_{0}M_{0}$, 
$$
\|\F_{\varepsilon,\rho}(\varphi)\|_{C^{0}}\le \|\F_{\varepsilon,\rho}(0)\|_{C^{0}}+\|\F_{\varepsilon,\rho}(\varphi)-\F_{\varepsilon,\rho}(0)\|_{C^{0}}\le M_{0}|\varepsilon|^{\widehat\gamma}+\frac{1}{8C_{0}}\|\varphi\|_{C^{2}}\le\frac{\mu_{0}|\varepsilon|^{\widehat\gamma}}{4C_{0}}
$$
and then, by \eqref{M1}, $\|\mathcal{Q}_{\varepsilon,\rho}(\varphi)\|_{C^{2}}\le \mu_{0}|\varepsilon|^{\widehat\gamma}$, namely, $\mathcal{Q}_{\varepsilon,\rho}(\mathcal{N}_{\varepsilon})\subset\mathcal{N}_{\varepsilon}$. Hence the assumptions of the contraction principle are satisfied and we can conclude that $\mathcal{Q}_{\varepsilon,\rho}$ admits a fixed point in $\mathcal{N}_{\varepsilon}$. 
\QED

\section{The variational argument}\label{var}

In this Section we complete the proof of Theorem \ref{mainteo}. The starting point is the result stated in Theorem \ref{T:reduced}, according to which for every $|\varepsilon|\ne 0$ small enough and for every $\rho \in [a_{1}|\varepsilon|^{-\delta},a_{2}|\varepsilon|^{-\delta}]$ there exists $\varphi_{\varepsilon,\rho} \in C^{2}(\R/2\pi\mathbb{Z})$ with $\|\varphi_{\varepsilon,\rho}\|_{C^2}\le\mu_{0}|\varepsilon|^{\min\{\gamma\delta,1-\delta\}}$ satisfying
$$
\K_{\varepsilon,\rho}(\varphi_{\varepsilon,\rho}) - B(u_{\varepsilon,\rho} + \varphi_{\varepsilon,\rho}\n_{\varepsilon,\rho}) = \lambda^1_{\varepsilon,\rho}\cos t + \lambda^2_{\varepsilon,\rho}\sin t
$$
where the Lagrange multipliers $\lambda^i_{\varepsilon,\rho} = \lambda^i_{\varepsilon,\rho}(\varphi_{\varepsilon,\rho})$ are given by \eqref{lambrappr}. Up to now, $a_1, a_2, \delta\in\R$ are arbitrary fixed constants with $0<a_1\leq a_2$ and $\delta \in (0,1)$, and $\mu_{0}$ is a constant depending on $a_1$, $a_2$ and $\delta$ but not on $\varepsilon$ neither on $\rho$. 

Here we prove that for a suitable choice of $a_1, a_2, \delta\in\R$ with $0<a_1< a_2$ and $\delta \in (0,1)$, for every $\varepsilon\ne 0$ small enough and with a suitable sign (the same of the coefficient $A$ in the assumption $(B_{1})$) one can find $\rho_{\varepsilon}\in [a_{1}|\varepsilon|^{-\delta},a_{2}|\varepsilon|^{-\delta}]$ such that $\lambda^1_{\varepsilon,\rho}=\lambda^2_{\varepsilon,\rho}=0$ when $\rho=\rho_{\varepsilon}$. To this aim, we exploit the variational nature of equation (\ref{eq:curvature}) which in fact corresponds to the Euler-Lagrange equation associated to a certain energy functional. 

The present Section consists in three parts: firstly we introduce the energy functional associated to (\ref{eq:curvature}) and we discuss its properties useful for the sequel. Then we tackle the equation $\lambda^1_{\varepsilon,\rho}=0$ and we prove that is admits a solution $\rho=\rho_{\varepsilon}$ under some conditions. Finally we show that as soon as $\lambda^1_{\varepsilon,\rho}=0$, then also $\lambda^2_{\varepsilon,\rho}=0$.

\subsection{The energy functional}

Since $B$ is radially symmetric of class $C^{1}$, there exists a $C^{1}$ map $b\colon[0,\infty)\to\R$ such that $B(v)=b(|v|)$ for every $v\in\R^{2}$. Set
$$
Q(v)=\frac{v}{|v|^{2}}\int_{0}^{|v|}b(s)s\de s\quad(v\in\R^{2}\setminus\{0\})~\!.
$$
Then $Q$ is a vector field with a continuous extension on $\R^{2}$, of class $C^{1}$ in $\R^{2}\setminus\{0\}$ and satisfying
\begin{gather}
\nonumber
\mathrm{div}~\!Q=B\text{~~on~$\R^{2}$}\\
\label{eq:Q-rotation}
Q(e^{i\theta}v) = e^{i\theta}Q(v) \quad \forall\theta\in\R~\!,~~ \forall v\in\R^2.
\end{gather}
For every $\varepsilon\in\R$, $\varepsilon\ne 1$, let 
\begin{gather*}
W^{1,1}_{\varepsilon}:=\{u\in W^{1,1}_{loc}(\R,\R^{2})~|~e^{\frac{2\pi i}{1-\varepsilon}}u(t)=u(t+2\pi)~\forall t\in\R~\!\}\\
\Omega_{\varepsilon}:=\{u\in W^{1,1}_{\varepsilon}~|~u(t)\ne 0~\forall t\in\R~\!,~\dot u\ne 0~\text{ a.e. in }\R~\!\}~\!.
\end{gather*}
Notice that $u_{\varepsilon,\rho}+\varphi\n_{\varepsilon,\rho}\in \Omega_{\varepsilon}$, as $\rho\in S_{\varepsilon}$, $\varphi\in \mathcal{N}_{\varepsilon}$ and $|\varepsilon|$ small enough. Set
$$
\mathscr{E}(u):=\int_{0}^{2\pi}\left(|\dot u|+Q(u)\cdot i\dot u\right)\de t\quad\forall u\in\Omega_{\varepsilon}~\!.
$$
One has that:
\begin{Lemma}
\label{L:energy}
The functional $\mathscr{E}$ is of class $C^{1}$ in $\Omega_{\varepsilon}$ and
$$
\mathscr{E}'(u)[h]=\int_{0}^{2\pi}\bigg(\frac{\dot u\cdot\dot h}{|\dot u|}+B(u)i\dot u\cdot h\bigg)\de t\quad\forall u\in\Omega_{\varepsilon}~\!,~\forall h\in W^{1,1}_{\varepsilon}.
$$
If in addition $u\in C^{2}$, then 
$$
\mathscr{E}'(u)[h]=\int_{0}^{2\pi}\left(B(u)-\mathcal{K}(u)\right)i\dot u\cdot h\de t\quad\forall h\in W^{1,1}_{\varepsilon}
$$
with $\mathcal{K}$ as in \eqref{mathcal-K}. 
\end{Lemma}

\Proof We can write 
$$
\mathscr{E}(u)=\mathscr{L}(u)+\mathscr{A}(u)\quad\text{where}\quad \mathscr{L}(u):=\int_{0}^{2\pi}|\dot u|\de t\quad\text{and}\quad\mathscr{A}(u):=\int_{0}^{2\pi}Q(u)\cdot i\dot u\de t~\!.
$$
In a standard way one shows that $\mathscr{L}$ is of class $C^{1}$ in $\Omega_{\varepsilon}$, and 
$$
\mathscr{L}'(u)[h]=\int_{0}^{2\pi}\frac{\dot u\cdot\dot h}{|\dot u|}\de t\quad\forall u\in\Omega_{\varepsilon}~\!,~\forall h\in W^{1,1}_{\varepsilon}~\!.
$$
Moreover, if $u\in C^{2}$, an integration by parts yields
$$
\mathscr{L}'(u)[h]=\frac{\dot u(2\pi)\cdot h(2\pi)}{|\dot u(2\pi)|}-\frac{\dot u(0)\cdot h(0)}{|\dot u(0)|}-\int_{0}^{2\pi}\left(\frac{\ddot u}{|\dot u|}-\frac{(\dot u\cdot\ddot u)\dot u}{|\dot u|^{3}}\right)\cdot h\de t~\!.
$$
Since $\dot u(2\pi)=e^{\frac{2\pi i}{1-\varepsilon}}\dot u(0)$ and $h(2\pi)=e^{\frac{2\pi i}{1-\varepsilon}}h(0)$, and using the decomposition 
$$
h=\frac{\dot u\cdot h}{|\dot u|^{2}}\dot u+\frac{i\dot u\cdot h}{|\dot u|^{2}}i\dot u
$$
one obtains
$$
\mathscr{L}'(u)[h]=-\int_{0}^{2\pi}\frac{(\ddot u\cdot i\dot u)(i\dot u\cdot h)}{|\dot u|^{3}}\de t\quad\forall h\in W^{1,1}_{\varepsilon}~\!.
$$
Let us study the regularity of the functional $\mathscr{A}$. Fix $u\in\Omega_{\varepsilon}$ and $h\in W^{1,1}_{\varepsilon}$. For $\epsilon\ne 0$ small enough one has that $u+\epsilon h\in\Omega_{\varepsilon}$ and there exists $r_{0}>0$ such that $|u(t)+\epsilon h(t)|> r_{0}>0$ for all $t\in\R$. One can write
$$
\frac{\mathscr{A}(u+\epsilon h)-\mathscr{A}(u)}{\epsilon}=\int_{0}^{2\pi}\frac{Q(u+\epsilon h)-Q(u)}{\epsilon}\cdot i\dot u\de t+\int_{0}^{2\pi}Q(u+\epsilon h)\cdot i\dot h\de t~\!.
$$
Since $Q\in C^{1}(\R^{2}\setminus\{0\},\R^{2})$ and, in particular, $Q$ and $DQ$ are locally uniformly continuous in $\R^{2}\setminus\{0\}$, by standard arguments, using also the embedding of $W^{1,1}(0,2\pi)$ into $C([0,2\pi])$, one can prove that
$$
\lim_{\epsilon\to 0}\frac{Q(u+\epsilon h)-Q(u)}{\epsilon}=DQ(u)h\quad\text{and}\quad \lim_{\epsilon\to 0}Q(u+\epsilon h)=Q(u)\quad\text{uniformly in }[0,2\pi].
$$
Therefore, by the Lebesgue Dominated Convergence Theorem, we infer that
\begin{equation*}
\begin{split}
\lim_{\epsilon\to 0}&\frac{\mathscr{A}(u+\epsilon h)-\mathscr{A}(u)}{\epsilon}=\int_{0}^{2\pi}DQ(u)[h,i\dot u]\de t+\int_{0}^{2\pi}Q(u)\cdot i\dot h\de t\\
&=\int_{0}^{2\pi}DQ(u)[h,i\dot u]\de t-\int_{0}^{2\pi}DQ(u)[\dot u,ih]\de t+Q(u(2\pi))\cdot ih(2\pi)-Q(u(0))\cdot ih(0)\\
&=\int_{0}^{2\pi}B(u)h\cdot i\dot u\de t~\!,
\end{split}
\end{equation*}
where the second equality is obtained by integration by parts, whereas, for the last one, one uses $u(2\pi)=e^{\frac{2\pi i}{1-\varepsilon}}u(0)$, $h(2\pi)=e^{\frac{2\pi i}{1-\varepsilon}}h(0)$, (\ref{eq:Q-rotation}) and the algebraic identity
$$
Mv\cdot iw-Mw\cdot iv=(\mathrm{tr}~\!M)v\cdot iw\quad\forall v,w\in\R^{2},
$$
where $M$ is any $2\times 2$ matrix and $\mathrm{tr}~\!M$ denotes 
its trace. The remaining part of the result can be proved in a standard way, following the same procedure as in \cite{CalGui} (see also \cite{CorPhD}). 
\QED

\subsection{The equation $\lambda^1_{\varepsilon,\rho}=0$}

Fix $a_{1},a_{2}>0$ such that
\begin{equation}
\label{eq:a1-a2}
0<a_{1}<\left(\frac{|A|\gamma}{2}\right)^{\frac{1}{\gamma+2}}<a_{2}
\end{equation}
and 
$$
\delta=\frac{1}{\gamma+2}
$$
and let $\varepsilon_{0}$ 
 be given by Theorem \ref{T:reduced}.
One has:

\begin{Lemma}\label{L:lambda1=zero}
There exists $\overline\varepsilon\in(0,\varepsilon_{0}]$ such that if $A>0$ then for every $\varepsilon\in(0,\overline\varepsilon)$ there exists $\rho_\varepsilon\in S_{\varepsilon} \big[a_{1}|\varepsilon|^{-\frac{1}{\gamma+2}},a_{2}|\varepsilon|^{-\frac{1}{\gamma+2}}\big]$ for which $\lambda^1_{\varepsilon,\rho_{\varepsilon}}=0$. If $A<0$ the same conclusion holds for every $\varepsilon\in(-\overline\varepsilon,0)$.
\end{Lemma}

\Proof
According to \eqref{lambrappr}, the equation $\lambda^1_{\varepsilon,\rho}=0$ can be written in the form
\begin{equation}\label{Kvarest}
\int_0^{2\pi}\K_{\varepsilon,\rho}(\varphi_{\varepsilon,\rho})\cos t \de t = \int_0^{2\pi}B(u_{\varepsilon,\rho} + \varphi_{\varepsilon,\rho} \n_{\varepsilon,\rho})\cos t \de t~\!.
\end{equation}
Our goal is to show that, taking $\delta=\frac{1}{\gamma+2}$, one has:
\begin{gather}
\label{lhs}
\int_0^{2\pi}\K_{\varepsilon,\rho}(\varphi_{\varepsilon,\rho})\cos t \de t =-2\pi\varepsilon\rho+F_{1}(\varepsilon,\rho)\\
\label{rhs}
\int_0^{2\pi}B(u_{\varepsilon,\rho} + \varphi_{\varepsilon,\rho} \n_{\varepsilon,\rho})\cos t \de t=-A\gamma\pi\rho^{-\gamma-1}+F_{2}(\varepsilon,\rho)
\end{gather}
with $F_{i}$ continuous functions such that
\begin{equation}
\label{lrhs}
\sup_{\rho\in S_{\varepsilon}}|F_{i}(\varepsilon,\rho)|=o\left(|\varepsilon|^{1-\delta}\right)\quad\text{as~~}\varepsilon\to 0~\!,~~i=1,2~\!.
\end{equation}
Assuming for a moment that \eqref{lhs}--\eqref{lrhs} hold true, let us complete the proof of the lemma. From \eqref{lhs}--\eqref{lrhs}, dividing by $\varepsilon\rho$, equation \eqref{Kvarest} becomes
\begin{equation}
\label{reduced-equation}
2-\frac{A\gamma}{\varepsilon\rho^{\gamma+2}}=F(\varepsilon,\rho)
\end{equation}
with $F$ continuous function such that 
\begin{equation}
\label{Fzero}
\sup_{\rho\in S_{\varepsilon}}|F(\varepsilon,\rho)|\to 0\quad\text{as~~}\varepsilon\to 0~\!.
\end{equation} 
Observe that for $\rho\in S_{\varepsilon}$ one has that $|\varepsilon|\rho^{\gamma+2}\le a_{2}^{\gamma+2}$. 
Hence, in order that \eqref{reduced-equation} admits a solution, $\varepsilon$ must have the same sign as $A$. Considering the case $A>0$ and defining $G_{\varepsilon}(\rho):=2-\frac{A\gamma}{\varepsilon\rho^{\gamma+2}}$, by \eqref{eq:a1-a2} we have that
$$
G_{\varepsilon}(a_{1}|\varepsilon|^{-\delta})=2-A\gamma a_{1}^{-(\gamma+2)}=:\alpha_{-}<0<\alpha_{+}:=2-A\gamma a_{2}^{-(\gamma+2)}=G_{\varepsilon}(a_{2}|\varepsilon|^{-\delta})~\!.
$$
Since $\alpha_{\pm}$ are independent of $\varepsilon$ and \eqref{Fzero} holds, there exists $\overline\varepsilon\in(0,\varepsilon_{0}]$ such that 
\begin{equation}
\label{bdry-sign}
G_{\varepsilon}(a_{1}|\varepsilon|^{-\delta})-F(\varepsilon,a_{1}|\varepsilon|^{-\delta})<0<G_{\varepsilon}(a_{2}|\varepsilon|^{-\delta})-F(\varepsilon,a_{2}|\varepsilon|^{-\delta})\quad\forall \varepsilon\in(0,\overline\varepsilon)~\!.
\end{equation}
Since the mapping $\rho\mapsto G_{\varepsilon}(\rho)-F(\varepsilon,\rho)$ is continuous on $S_{\varepsilon}$, by \eqref{bdry-sign}, it must vanish at some $\rho_{\varepsilon}\in S_{\varepsilon}$, for every $\varepsilon\in(0,\overline\varepsilon)$, namely \eqref{reduced-equation}, and thus also \eqref{Kvarest} are satisfied when $\rho=\rho_{\varepsilon}$. If $A<0$ one can repeat the same argument, taking $\varepsilon\in(-\overline\varepsilon,0)$. Thus the lemma is proved. It remains to check \eqref{lhs} and \eqref{rhs}. 
\medskip

\noindent\textit{\textbf{Proof of \eqref{lhs}.}}
Setting 
\begin{equation}\label{resto}
\mathcal{R}_{\varepsilon,\rho}(\varphi):=\K_{\varepsilon,\rho}(\varphi)-\K_{\varepsilon,\rho}(0)-\K'_{\varepsilon,\rho}(0)[\varphi]\quad\forall\varphi\in\mathcal{N}~\!,
\end{equation}
we can write
\begin{equation}\label{lhs1}
\begin{aligned}
&\int_0^{2\pi}\K_{\varepsilon,\rho}(\varphi_{\varepsilon,\rho})\cos t \de t=\int_0^{2\pi}(\mathcal{L}_{0}\varphi_{\varepsilon,\rho})\cos t\de t+\int_0^{2\pi}(\K'_{\varepsilon,\rho}(0)-\mathcal{L}_{0})[\varphi_{\varepsilon,\rho}]\cos t\de t\\
&\qquad+\int_0^{2\pi}\mathcal{R}_{\varepsilon,\rho}(\varphi_{\varepsilon,\rho})\cos t\de t
-\varepsilon\rho\int_0^{2\pi}(\K_{\varepsilon,\rho}(0)-1)\de t
+\int_0^{2\pi}(\K_{\varepsilon,\rho}(0)-1)(\varepsilon\rho+\cos t)\de t~\!.
\end{aligned}
\end{equation}
Integrating by parts twice we readily get that
\begin{equation}\label{lhs2}
\int_0^{2\pi} \left(\Le_{0} \varphi_{\varepsilon,\rho} \right)\cos t \de t = \int_0^{2\pi}\varphi_{\varepsilon,\rho}\left(\Le_{0} \cos t\right)\de t = 0~\!, 
\end{equation}
because $\cos t \in \ker(\Le_{0})$. Moreover, by \eqref{Lnest}, since $\rho\in S_{\varepsilon}$ and $\varphi_{\varepsilon,\rho}$ satisfies \eqref{estimphi}, we have
$$
\left|\int_0^{2\pi}(\K'_{\varepsilon,\rho}(0)-\mathcal{L}_{0})[\varphi_{\varepsilon,\rho}]\cos t\de t\right|\le C|\varepsilon|^{1-\delta+\widehat\gamma}~\!.
$$
Therefore we can write
\begin{equation}\label{lhs3}
\int_0^{2\pi}(\K'_{\varepsilon,\rho}(0)-\mathcal{L}_{0})[\varphi_{\varepsilon,\rho}]\cos t\de t=|\varepsilon|^{1-\delta+\widehat\gamma}K_{1}(\varepsilon,\rho)
\end{equation}
where $K_{1}$ is a continuous function, uniformly bounded with respect to $\rho\in S_{\varepsilon}$ and $\varepsilon\in[-\varepsilon_{0},\varepsilon_{0}]\setminus\{0\}$. 
In order to estimate the third integral on the right-hand side of \eqref{lhs1}, we observe that, fixing $\varphi$ and $t\in\R$, by \eqref{resto} we can write
$$
\mathcal{R}_{\varepsilon,\rho}(\varphi)(t)=\mathcal{R}_{\varepsilon,\rho}(0)(t)+\int_{0}^{1}\mathcal{R}'_{\varepsilon,\rho}(s\varphi)[\varphi](t)\de s=\int_{0}^{1}\left(\mathcal{K}'_{\varepsilon,\rho}(s\varphi)-\mathcal{K}'_{\varepsilon,\rho}(0)\right)[\varphi](t)\de s~\!.
$$
Hence, taking $\varphi=\varphi_{\varepsilon,\rho}$ and using \eqref{estimphi} and \eqref{abc0}, we obtain
$$
\|\mathcal{R}_{\varepsilon,\rho}(\varphi_{\varepsilon,\rho})\|_{C^{0}}\le\sup_{s\in[0,1]}\left\|\mathcal{K}'_{\varepsilon,\rho}(s\varphi_{\varepsilon,\rho})-\mathcal{K}'_{\varepsilon,\rho}(0)\right\|\|\varphi_{\varepsilon,\rho}\|_{C^{2}}\le C|\varepsilon|^{2\widehat\gamma}~\!.
$$
This estimate allows us to write
\begin{equation}\label{lhs4}
\int_0^{2\pi}\mathcal{R}_{\varepsilon,\rho}(\varphi_{\varepsilon,\rho})\cos t\de t=|\varepsilon|^{2\widehat\gamma}K_{2}(\varepsilon,\rho)
\end{equation}
where $K_{2}$ is a continuous function, uniformly bounded with respect to $\rho\in S_{\varepsilon}$ and $\varepsilon\in[-\varepsilon_{0},\varepsilon_{0}]\setminus\{0\}$. By \eqref{passag1948} we have that
$$
\left|\int_0^{2\pi}(\K_{\varepsilon,\rho}(0)-1)\de t\right|\le C|\varepsilon|^{1-\delta}
$$
and then
\begin{equation}\label{lhs5}
\varepsilon\rho\int_0^{2\pi}(1-\K_{\varepsilon,\rho}(0))\de t=|\varepsilon|^{2(1-\delta)}K_{3}(\varepsilon,\rho)
\end{equation}
where $K_{3}$ is a continuous function, uniformly bounded with respect to $\rho\in S_{\varepsilon}$ and $\varepsilon\in[-\varepsilon_{0},\varepsilon_{0}]\setminus\{0\}$. In order to estimate the last integral on the right-hand side of \eqref{lhs1} we apply Lemma \ref{L:energy}, with $B\equiv 1$ and $Q(v)=\frac{v}{2}$. In particular we set
$$
\E_{0}(u):=\mathscr{L}(u)+\mathscr{A}_{0}(u)\quad\text{where}~~\mathscr{A}_{0}(u)=\frac{1}{2}\int_{0}^{2\pi}u\cdot i\dot u\de t~\!.
$$
By Lemma \ref{L:energy} we have that 
$$
\frac{\partial}{\partial\rho}\left[\E_{0}(u_{\varepsilon,\rho})\right]=\E'_{0}(u_{\varepsilon,\rho})[h_{\varepsilon}]
$$
where $h_{\varepsilon}(t)=e^{\frac{i\varepsilon t}{1-\varepsilon}}$. Since $u_{\varepsilon,\rho}, h_{\varepsilon}\in\Omega_{\varepsilon}$, we obtain
\begin{equation}\label{radiusderiv}
(1-\varepsilon)\frac{\partial}{\partial\rho}\left[\E_{0}(u_{\varepsilon,\rho})\right]=\int_{0}^{2\pi}[\mathcal{K}(u_{\varepsilon,\rho})-1](\varepsilon\rho+\cos t)\de t
\end{equation}
with $\mathcal{K}$ as in \eqref{mathcal-K}.
We claim that 
\begin{equation}\label{est3}
(1-\varepsilon)\frac{\partial}{\partial\rho}\left[\E_{0}(u_{\varepsilon,\rho})\right]= - 2\pi \varepsilon \rho+ |\varepsilon|^{2-\delta}K_4(\varepsilon,\rho)\quad\forall\rho\in S_{\varepsilon}~\!,\quad\forall\varepsilon\in[-\varepsilon_{0},\varepsilon_{0}]\setminus\{0\}~\!.
\end{equation}
where $K_{4}(\varepsilon,\rho)$ is a continuous function, uniformly bounded with respect to $\rho\in S_{\varepsilon}$ and $\varepsilon\in[-\varepsilon_{0},\varepsilon_{0}]\setminus\{0\}$. Therefore, by \eqref{lhs1}--\eqref{est3}, since $\min\{1-\delta+\widehat\gamma,~\!2\widehat\gamma,~\!2-2\delta,~\!2-\delta\}>1-\delta$ for $\delta=\frac{1}{\gamma+2}$ and $\gamma>1$, \eqref{lhs} is true. 
\medskip

\noindent
\emph{Proof of (\ref{est3}).}
Recall that $\E_{0}(u):=\mathscr{L}(u)+\mathscr{A}_{0}(u)$. 
A simple computation shows that
\begin{equation}
\label{eq:K31}
(1-\varepsilon)\frac{\partial}{\partial\rho}[\A_{0}(u_{\varepsilon,\rho})] = - 2\pi \varepsilon \rho~\!.
\end{equation}
It remains to estimate the part of the energy associated to the length functional $\mathscr{L}$. A direct computation yields
$$
(1-\varepsilon)\mathscr{L}(u_{\varepsilon,\rho})=\int_0^{2\pi}\sqrt{1+\varepsilon^2\rho^2+2\varepsilon\rho\cos t}\de t~\!.
$$
Notice that $\varepsilon\in[-\varepsilon_{0},\varepsilon_{0}]\setminus\{0\}$ and $\rho \in S_\varepsilon$ implies $|\varepsilon \rho|\le a_{2}\varepsilon_{0}^{1-\delta}=:r_{0}<1$, thanks to \eqref{nsotto}, and the integral on the right-hand side in the above formula is well defined. It is convenient to introduce the auxiliary function
$$
G(r):=\int_0^{2\pi}\sqrt{1+2r\cos t+r^{2}}\de t\quad\forall|r|<1~\!.
$$
A threefold application of the derivation theorem for integrals depending on a parameter yields that $G$ is of class $C^{3}$ in $(-1,1)$, with
$$
G^{(k)}(r)=\int_0^{2\pi}\frac{\partial^{k}g}{\partial r^{k}}(r,t)\de t\quad\forall|r|<1~~(k=1,2,3)\text{,~~where~~}g(r,t)=\sqrt{1+2r\cos t+r^{2}}~\!.
$$
In particular $G'(0)=0$, $G''(0)=\pi$ and $\frac{\partial^{3}g}{\partial r^{3}}$ is uniformly bounded with respect to $t\in[0,2\pi]$ and $|r|\le r_{0}$. Therefore
$$
G'(r)=\pi r+r^{2}H(r)\quad\forall r\in[-r_{0},r_{0}]
$$
where $H$ is a continuous function on $[-r_{0},r_{0}]$. Since $(1-\varepsilon)\mathscr{L}(u_{\varepsilon,\rho})=G(\varepsilon\rho)$, $\rho\in S_{\varepsilon}$ and $\varepsilon\in[-\varepsilon_{0},\varepsilon_{0}]\setminus\{0\}$, we obtain that
\begin{equation}
\label{eq:K32}
(1-\varepsilon)\frac{\partial}{\partial\rho}[\mathscr{L}(u_{\varepsilon,\rho})]=\pi\varepsilon\left[\varepsilon\rho+|\varepsilon|^{2(1-\delta)}\widetilde{H}(\varepsilon,\rho)\right]\quad\forall\varepsilon\in[-\varepsilon_{0},\varepsilon_{0}]\setminus\{0\}~\!,\quad\forall\rho\in S_{\varepsilon}~\!,
\end{equation}
where $\widetilde{H}(\varepsilon,\rho)=\rho^{2}|\varepsilon|^{2\delta}H(\varepsilon\delta)$ is a continuous function, bounded 
uniformly with respect to $\varepsilon\in[-\varepsilon_{0},\varepsilon_{0}]\setminus\{0\}$ and $\rho\in S_{\varepsilon}$. Thus (\ref{est3}) follows from (\ref{eq:K31}) and (\ref{eq:K32}).
\medskip

\noindent\textit{\textbf{Proof of \eqref{rhs}.}} It is convenient to write 
\begin{equation}
\label{eq:B-decomposition}
B(v)=1+\frac{A}{|v|^{\gamma}}+\frac{B_{1}(v)}{|v|^{\gamma + \beta}} \end{equation}
where
$$
\beta:=\min\{1,\gamma_{1}-\gamma\}\in(0,1]
$$
and $B_{1}\colon\R^{2}\to\R$ is a continuous radially symmetric function, defined by (\ref{eq:B-decomposition}), satisfying
$$
B_{1}(v)=\left\{\begin{array}{ll}o(1)&\text{if $\gamma_{1}>\gamma+1$}\\
A_{1}+o(|v|^{\beta-1})&\text{if $\gamma_{1}\le\gamma+1$}\end{array}\right.\text{as $|v|\to\infty$.}
$$
By \eqref{eq:B-decomposition}, we have that
\begin{equation}\label{Hdeco}
\begin{split}
\int_0^{2\pi}&B(u_{\varepsilon,\rho} + \varphi_{\varepsilon,\rho} \n_{\varepsilon,\rho})\cos t \de t=A\int_0^{2\pi}\frac{\cos t}{|u_{\varepsilon,\rho}|^\gamma}\de t\\
&+ A\int_0^{2\pi}\cos t \left(\frac{1}{|u_{\varepsilon,\rho} + \varphi_{\varepsilon,\rho}\n_{\varepsilon,\rho}|^\gamma} - \frac{1}{|u_{\varepsilon,\rho}|^\gamma}\right)\de t
+ \int_0^{2\pi}\frac{B_{1}(u_{\varepsilon,\rho} + \varphi_{\varepsilon,\rho}\n_{\varepsilon,\rho})}{|u_{\varepsilon,\rho} + \varphi_{\varepsilon,\rho}\n_{\varepsilon,\rho}|^{\gamma+\beta}}\cos t \de t~\!.
\end{split}
\end{equation}
In order to estimate the first integral on the right-hand side of \eqref{Hdeco}, we write
$$
\int_0^{2\pi}\frac{\cos t}{|u_{\varepsilon,\rho}|^\gamma}\de t=\rho^{-\gamma}G_{1}\big(\rho^{-1}\big)\quad\text{where~~}G_{1}(r)=\int_{0}^{2\pi}\frac{\cos t}{(1+2r\cos t+r^{2})^{\frac{\gamma}{2}}}\de t\!~.
$$
Notice that as $\rho\in S_{\varepsilon}$ and $\varepsilon\in(0,\varepsilon_{0}]$, one has that $0<\rho^{-1}\le a_{1}^{-1}\varepsilon_{0}^{\delta}=:r_{1}<1$. 
A double application of the derivation theorem for integrals depending on a parameter yields that $G_{1}$ is of class $C^{2}$ in $(-1,1)$, with
$$
G'_{1}(r)=\int_{0}^{2\pi}\frac{\partial g_{1}}{\partial r}(r,t)\de t\quad\text{and}\quad G''_{1}(r)=\int_{0}^{2\pi}\frac{\partial^{2} g_{1}}{\partial r^{2}}(r,t)\de t\quad\text{where}\quad g_{1}(r,t)=\frac{\cos t}{(1+2r\cos t+r^{2})^{\frac{\gamma}{2}}}~\!.
$$
In particular $G_{1}(0)=0$, $G'_{1}(0)=-\gamma\pi$ and $\frac{\partial^{2} g_{1}}{\partial r^{2}}$ is uniformly bounded with respect to $t\in[0,2\pi]$ and $|r|\le r_{1}$. Therefore
$$
G_{1}(r)=-\gamma\pi r+r^{2}H_{1}(r)\quad\forall r\in[-r_{1},r_{1}]
$$
where $H_{1}$ is a continuous function on $[-r_{1},r_{1}]$. As a consequence
\begin{equation}
\label{H1dec}
\int_0^{2\pi}\frac{\cos t}{|u_{\varepsilon,\rho}|^\gamma}\de t 
=-\gamma\pi\rho^{-\gamma-1}+\rho^{-\gamma-2}H_{1}(\rho^{-1})=-\gamma\pi\rho^{-\gamma-1}+|\varepsilon|^{\delta(\gamma+2)}
\widetilde{B}_{1}(\varepsilon,\rho)
\end{equation}
where $\widetilde{B}_{1}(\rho)$ is a continuous function, bounded 
uniformly with respect to $\varepsilon\in[-\varepsilon_{0},\varepsilon_{0}]\setminus\{0\}$ and $\rho\in S_{\varepsilon}$.
The second integral on the right-hand side of \eqref{Hdeco} can be estimated as follows. By the mean value theorem, and using \eqref{estimphi} and \eqref{passag17}, for every $t$ there exists $s\in[0,1]$ such that 
$$
\big||u_{\varepsilon,\rho}(t) + \varphi_{\varepsilon,\rho}(t)\n_{\varepsilon,\rho}(t)|^{-\gamma} - |u_{\varepsilon,\rho}(t)|^{-\gamma}\big|
\leq \gamma |u_{\varepsilon,\rho}(t) + s \varphi_{\varepsilon,\rho}(t)\n_{\varepsilon,\rho}(t)|^{-\gamma-1}\|\varphi_{\varepsilon,\rho}\|_{C^2} \le C|\varepsilon|^{\delta(\gamma +1)+ \widehat\gamma}.
$$
Therefore
\begin{equation}
\label{H2dec}
\int_0^{2\pi}\left(\frac{1}{|u_{\varepsilon,\rho} + \varphi_{\varepsilon,\rho}\n_{\varepsilon,\rho}|^\gamma} - \frac{1}{|u_{\varepsilon,\rho}|^\gamma}\right)\cos t \de t=|\varepsilon|^{\delta(\gamma+1) +\widehat\gamma}\widetilde{B}_2(\varepsilon,\rho)
\end{equation}
where $\widetilde{B}_2$ is a continuous function, uniformly bounded with respect to $\varepsilon\in[-\varepsilon_{0},\varepsilon_{0}]\setminus\{0\}$ and $\rho\in S_{\varepsilon}$. 
The estimate of the last integral on the right-hand side of \eqref{Hdeco} can be accomplished in different ways according that $\gamma_{1}>\gamma+1$ or $\gamma_{1}\le\gamma+1$. Let us examine firstly the case $\gamma_{1}>\gamma+1$, in which $\beta=1$ and $B_{1}(v)=o(|v|)$ as $|v|\to\infty$. 
Since $|u_{\varepsilon,\rho}(t) + \varphi_{\varepsilon,\rho}(t)\n_{\varepsilon,\rho}(t)|\ge C|\varepsilon|^{-\delta}$ for every $t$, we have that
$$
\left|\int_0^{2\pi}\frac{B_{1}(u_{\varepsilon,\rho} + \varphi_{\varepsilon,\rho}\n_{\varepsilon,\rho})}{|u_{\varepsilon,\rho} + \varphi_{\varepsilon,\rho}\n_{\varepsilon,\rho}|^{\gamma+\beta}}\cos t \de t\right|\le C|\varepsilon|^{\delta(\gamma+1)}\sup_{|v|\ge C|\varepsilon|^{-\delta}}|B_{1}(v)|~\!.
$$
For $\delta=\frac{1}{\gamma+2}$, one has that $\delta(\gamma+1)=1-\delta$ and 
\begin{equation}
\label{H3dec}
\int_0^{2\pi}\frac{B_{1}(u_{\varepsilon,\rho} + \varphi_{\varepsilon,\rho}\n_{\varepsilon,\rho})}{|u_{\varepsilon,\rho} + \varphi_{\varepsilon,\rho}\n_{\varepsilon,\rho}|^{\gamma+\beta}}\cos t \de t=|\varepsilon|^{1-\delta}\widetilde{B}_3(\varepsilon,\rho)
\end{equation}
where $\widetilde{B}_3$ is a continuous function, such that $\sup_{\rho\in S_{\varepsilon}}|\widetilde{B}_3(\varepsilon,\rho)|\to 0$ as $\varepsilon\to 0$. Therefore in this case, by \eqref{Hdeco}--\eqref{H3dec}, since $\min\{\delta(\gamma+2),~\!\delta(\gamma+1) +\widehat\gamma\}>1-\delta$ for $\delta=\frac{1}{\gamma+2}$ and $\gamma>1$, \eqref{rhs} is proved. Finally let us consider the case $\gamma_{1}\le\gamma+1$, in which $\beta\in(0,1]$ and $B_{1}(v)=A_{1}+|v|^{\beta-1}B_{2}(v)$ with $B_{2}(v)\to 0$ as $|v|\to\infty$. We split the last integral on the right-hand side of \eqref{Hdeco} as follows:
\begin{equation}
\label{second-case}
\begin{split}
\int_0^{2\pi}\frac{B_{1}(u_{\varepsilon,\rho} + \varphi_{\varepsilon,\rho}\n_{\varepsilon,\rho})}{|u_{\varepsilon,\rho} + \varphi_{\varepsilon,\rho}\n_{\varepsilon,\rho}|^{\gamma+\beta}}\cos t &\de t=A_{1}\int_0^{2\pi}\left[\frac{1}{|u_{\varepsilon,\rho} + \varphi_{\varepsilon,\rho}\n_{\varepsilon,\rho}|^{\gamma+\beta}} - \frac{1}{|u_{\varepsilon,\rho}|^{\gamma+\beta}}\right]\cos t \de t\\
&+ A_{1} \int_0^{2\pi}  \frac{\cos t}{|u_{\varepsilon,\rho}|^{\gamma+\beta}}\de t+ \int_0^{2\pi}  \frac{B_{2}(u_{\varepsilon,\rho} + \varphi_{\varepsilon,\rho}\n_{\varepsilon,\rho})}{|u_{\varepsilon,\rho} + \varphi_{\varepsilon,\rho}\n_{\varepsilon,\rho}|^{\gamma+1}}\cos t\de t~\!.
\end{split}
\end{equation}
The first two integrals on the right-hand side of \eqref{second-case} can be studied as \eqref{H2dec} and \eqref{H1dec}, respectively, with $\gamma+\beta$ instead of $\gamma$, whereas the third one is like \eqref{H3dec}, with $B_{2}$ instead of $B_{1}$. In conclusion, since $\gamma+\beta>\gamma$,  \eqref{H3dec} holds true even for $\gamma_{1}\le\gamma+1$. 
Hence, noting that $\min\{\delta(\gamma+1)+\widehat\gamma,~\!\delta(\gamma+2)\}>1-\delta$ for $\delta=\frac{1}{\gamma+2}$, \eqref{Hdeco}--\eqref{H3dec} imply \eqref{rhs}. 
\QED

\subsection{The equation $\lambda^2_{\varepsilon,\rho}=0$}

The fact that the second Lagrange multiplier $\lambda^2_{\varepsilon,\rho}$ vanishes whenever the first one does is discussed below and is a consequence of the invariance under rotation. 

\begin{Lemma}\label{L:lambda2=zero}
If $u=u_{\varepsilon,\rho}+\varphi\n_{\varepsilon,\rho}$ solves $\mathcal{K}(u)-B(u)=\lambda\sin t$, with $\mathcal{K}$ defined in \eqref{mathcal-K}, and $\varphi\in\mathcal{N}_{\varepsilon}$, then $\lambda=0$.
\end{Lemma}

\Proof 
Since $Q$ solves (\ref{eq:Q-rotation}) and $\big|\frac{\mathrm{d}}{\mathrm{d}t}(e^{i\theta}u)\big|=|\dot u|$ for every $\theta\in\R$, one has that
$$
\E(e^{i\theta}u)=\E(u)\quad\forall\theta\in\R~\!.
$$
Hence, by the $C^{1}$ regularity of $\E$ (Lemma \ref{L:energy}),
\begin{equation}
\label{eq:lambda2A}
0=\frac{\mathrm{d}}{\mathrm{d}\theta}\big[\E((e^{i\theta}u)\big]=\E'(e^{i\theta}u)[h]\quad\text{where~~}h=ie^{i\theta}\dot u~\!.
\end{equation}
As $e^{i\theta}u,h\in\Omega_{\varepsilon}$, by Lemma \ref{L:energy} and by the symmetry of $B$, one has that
\begin{equation}
\label{eq:lambda2B}
\E'(e^{i\theta}u)[h]=\int_{0}^{2\pi}[B(e^{i\theta}u)-\K(e^{i\theta}u)]ie^{i\theta}\dot u\cdot ie^{i\theta}u\de t=\int_{0}^{2\pi}[B(u)-\K(u)]\dot u\cdot u\de t~\!.
\end{equation}
Using the assumption $\K(u)-B(u)=\lambda\sin t$, integrating by parts and taking into account that $|u(2\pi)|=|u(0)|$, by \eqref{eq:lambda2A} and \eqref{eq:lambda2B} we obtain
\begin{equation}
\label{eq:lambda2C}
0=\lambda\int_{0}^{2\pi}|u|^{2}\cos t\de t~\!.
\end{equation}
Recalling that $u=u_{\varepsilon,\rho} + \varphi\n_{\varepsilon,\rho}$, with $|u_{\varepsilon,\rho}|^{2}=\rho^{2}+1+2\rho\cos t$ and $|\n_{\varepsilon,\rho}|=1$, from (\ref{eq:lambda2C}) it follows that
\begin{equation}
\label{eq:lambda2D}
\lambda\left[2\rho\int_{0}^{2\pi}\cos^{2}t\de t+\int_{0}^{2\pi}\left(\varphi^{2}+2\varphi u_{\varepsilon,\rho}\cdot\n_{\varepsilon,\rho}\right)\cos t\de t\right]=0~\!.
\end{equation}
Finally, we observe that, since $\varphi\in\mathcal{N}_{\varepsilon}$,  
$$
\left|\int_{0}^{2\pi}\left(\varphi^{2}+2\varphi u_{\varepsilon,\rho}\cdot\n_{\varepsilon,\rho}\right)\cos t\de t\right|\le C\left(|\varepsilon|^{2\min\{\gamma\delta,1-\delta\}}+(\rho+1)|\varepsilon|^{\min\{\gamma\delta,1-\delta\}}\right)<2\rho\int_{0}^{2\pi}\cos^{2}t\de t
$$
for $|\varepsilon|\ne 0$ small enough, and then (\ref{eq:lambda2D}) implies $\lambda=0$.
\QED
\bigskip

\noindent
\textbf{Acknowledgements.}
The authors are members of the Gruppo Nazionale per l'Analisi Matematica, la Probabilit\`a e le loro Applicazioni (GNAMPA) of the Istituto Nazionale di Alta Matematica (INdAM) and their research has been partly supported by the ERC Advanced Grant 2013 n. 339958 COMPAT -- Complex Patterns for Strongly Interacting Dynamical Systems. The research of the first author has been partly supported by the PRIN Project 2015KB9WPT ``Variational methods, with applications to problems in mathematical physics and geometry''.

\end{document}